\documentclass[moor]{informs1}              % for a regular run

\RequirePackage{fix-cm}
\usepackage{graphicx}
\usepackage{psfrag}
\usepackage{subfigure}
\usepackage{url}
\usepackage{color}
\usepackage{cite}
\usepackage{epsfig}
\usepackage{amsmath,amssymb}
\usepackage{slashbox}
\usepackage{longtable}

\newcommand{\cX}{{\cal X}}
\newcommand{\cY}{{\cal Y}}
\newcommand{\cZ}{{\cal Z}}
\newcommand{\eqdef}{\overset{\text{def}}{=}}
\newcommand{\R}{\mathbb{R}}

\newcommand{\Uconst}{{\color{black}\gamma}}
\newcommand{\Prob}{\mathbf{P}}
\newcommand{\Exp}[1]{\mathbf{E}\left[#1\right]}
\newcommand{\dist}{{\rm dist}}

\newcommand{\Ib}{\mathbb{I}}

\providecommand{\rset}[1]{\mathbb{R}^}

\providecommand{\norm}[1]{\lVert#1\rVert}

\usepackage[colorlinks=true,breaklinks=true,bookmarks=true,urlcolor=blue,
     citecolor=blue,linkcolor=blue,bookmarksopen=false,draft=false]{hyperref}

\def\EMAIL#1{\href{mailto:#1}{#1}}% When hyperref is used, otherwise outcomment
\TheoremsNumberedThrough     % Preferred (Theorem 1, Lemma 1, Theorem 2)
\EquationsNumberedThrough    % Default: (1), (2), ...
\def\TheoremsNumberedThrough{%
\theoremstyle{TH}%
\newtheorem{theorem}{Theorem}
\newtheorem{lemma}{Lemma}

\newtheorem{corollary}{Corollary}

\newtheorem{assumption}{Assumption}
\theoremstyle{EX}

\newtheorem{example}{Example}

\newtheorem{definition}{Definition}

}

\usecounter{algorithmenumi}

\begin{document}
\RUNTITLE{Necoara et al., Randomized projection methods for convex feasibility problems}

% Full title.
\TITLE{ Randomized projection methods for convex feasibility problems: conditioning and convergence rates}

\ARTICLEAUTHORS{
\AUTHOR{Ion Necoara and Andrei Patrascu}
\AFF{Automatic Control and
Systems Engineering Department, University Politehnica Bucharest,
060042 Bucharest, Romania, \EMAIL{ion.necoara@acse.pub.ro.}}

\AUTHOR{Peter Richtarik}
\AFF{School of Mathematics, The Maxwell Institute for
Mathematical Sciences, University of Edinburgh, United Kingdom,
\EMAIL{peter.richtarik@ed.ac.uk.} }

%\AUTHOR{Andrei P\u atra\c scu}
%\AFF{Automatic Control and
%Systems Engineering Department, University Politehnica Bucharest,
%060042 Bucharest, Romania, \EMAIL{andrei.patrascu@acse.pub.ro.} }
} 

\ABSTRACT{
Finding a point in the intersection of a collection of closed convex sets, that is the convex feasibility problem, represents the main modeling strategy for many computational problems. In this paper we analyze new stochastic reformulations of the  convex feasibility problem in order to facilitate the development of new  algorithmic schemes. We also analyze the conditioning problem parameters using certain (linear) regularity assumptions on the individual convex sets.
Then, we introduce a  general random projection algorithmic framework, which extends to the random settings many existing projection schemes, designed for the general convex feasibility problem.  Our general random projection algorithm allows to project simultaneously on several sets, thus providing great flexibility in matching the implementation of the algorithm on the parallel architecture at hand. Based on the conditioning parameters, besides the asymptotic convergence results, we also derive explicit sublinear and linear convergence rates for this general algorithmic framework.}

% \KEYWORDS{deterministic inventory theory; infinite linear programming duality;
%  existence of optimal policies; semi-Markov decision process; cyclic schedule}

% \MSCCLASS{Primary: 90B05; secondary: 90C40, 90C90}

% \ORMSCLASS{Primary: Inventory/production: deterministic multi-item;
%  secondary: dynamic programming/optimal control: deterministic
%  semi-Markov; programming: infinite dimensional}

\HISTORY{First version: March 2017.}

%\date{Received: March 2017 / Accepted: date}

\maketitle

%%%%%%%%%%%%%%%%%%%%%%%%%%%%%%%%%%%%%%%%%%%%%%%%%%%%%%%%%%%%%%%%

\section{Introduction}
\noindent Finding a point in the intersection of a collection of
closed convex sets, that is \textit{the convex feasibility problem},
represents a modeling paradigm which has been used for many decades
for posing and solving engineering and physics problems. Among the
most important applications modeled by the convex feasibility
formalism are: radiation therapy treatment planning \cite{HerChe:08},
computerized tomography \cite{Her:09} and magnetic
resonance imaging \cite{SamKho:04}; wavelet-based denoising
\cite{ChoBar:04}, color imaging \cite{Sha:00} and demosaicking
\cite{LuKar:10}; antenna design \cite{GuSta:04} and sensor networks
problems\cite{BlaHer:06}; data compression \cite{LieYan:05}, neural
networks \cite{StaYan:98} and adaptive filtering \cite{YukYam:06}.

\noindent  Convex feasibility problems have various formulations,
such as finding the fixed points of a nonexpansive operator, the set
of optimal solutions of a specific optimization problem or the set
of solutions to some convex inequalities.   Projection methods were
first used for solving systems of linear equalities  \cite{Kac:37}
or  linear inequalities \cite{MotSch:54}, and then extended to
general convex feasibility problems, e.g. in \cite{Com:96}.
Projection methods are very attractive in applications since they
are able to handle problems of huge dimensions with a very
large number of convex sets in the intersection. For instance, the \textit{projection algorithm} which represents one of the first iterative algorithms for feasibility problems, rely at each iteration on orthogonal projections onto given individual sets. Its simple algorithmic structure supports the current large scale setting and can be easily adapted to parallel environments,
making such schemes adequate to modern computational architectures.
If the iteration of a given projection algorithm rely on an alternating sequence
of projections onto sets over the iterations, then
it belongs to an \textit{alternating projection} schemes
\cite{BauBor:96, BauNol:13, GowRic:15, Ned:10}. Furthermore, depending on the
variant of the alternating projection algorithm, the current set (or
sets) on which the projection is made can be chosen, for example, in
a random, cyclic or greedy manner. Otherwise, if the scheme uses at
current iteration an average of multiple projections of  the current
iterate onto various sets, then also it can be viewed as an
\textit{average projection} algorithm \cite{CenElf:01, CenChe:12}.

\noindent The convergence properties, the iteration complexity and even the inherent limitations of the class of projection schemes has been intensely analyzed over
the last decades, as it can be seen in \cite{AusTeb:01, BecTeb:03, BecTeb:04, BecTeb:11, BauBor:96, BauNol:13,
 CenChe:12, CenElf:01, Com:96, GowRic:15, Ned:10, Ned:11} and the
references therein. In \cite{AusTeb:01} a barycenter type projection algorithm is developed, which allows the efficient handling of feasibility problems arising in the nonnegative orthant. The proposed method uses approximate projections on the sets and is proven globally convergent under the sole assumption that the given intersection is nonempty and the errors are controllable.
An important contribution is made in \cite{BecTeb:03}, where the rates of convergence of some projections algorithms are analyzed for solving the general convex feasibility problem. Besides revealing some connections between the Slater's condition and the classical linear regularity property,
the authors show that if the Slater's condition does not hold, the projection algorithms can behave quite badly, i.e. with a rate of convergence which
is not bounded. Moreover, the authors also propose an alternative local linear regularity bound  to derive further convergence rate results. Linear convergence of the conditional gradient method applied on the equivalent optimization formulation of the problem of finding a point in the intersection of an affine set with a compact convex set is derived in \cite{BecTeb:04}. In a more general setting, \cite{BecTeb:11} studies the problem of finding a point in the intersection of affine constraints with a nonconvex closed set and a simple gradient projection scheme is developed. The scheme is proven to converge to a
unique solution of the problem, at a linear rate, under a natural assumption defined in terms of the problem's data. 

\vspace{5pt}

\noindent \textit{Contributions}. Below, we clarify the relationship
and differences   between our work and earlier research in this
direction. In particular, the main contributions of this paper
consist in unifying and extending existing projection methods in
several aspects:

\noindent $(i)$  The classical convex feasibility problem was
usually formulated for a finite intersection of simple convex sets.
While  finding a point in the intersection of a finite number of
convex sets is a problem with its own challenges, it does not cover
many interesting applications modeled by an intersection of (infinite)
countable/uncountable number of simple convex sets (see e.g. \cite{PatNec:17}).  
In this paper we present several new equivalent stochastic formulations of the
convex feasibility problem, which allow us to deal with
intersections of families of convex sets that may be even
uncountable.

\noindent $(ii)$  From an algorithmic point of view, most of the
previous approaches are limited to cycle based alternating
projection schemes. Moreover, for this strategy it is difficult to
prove asymptotic convergence and to estimate the rate of convergence 
in the general convex feasibility case. Therefore, we introduce a  general random projection algorithmic framework, which covers or extends to the random
settings many existing projection schemes, designed for the general
convex feasibility problem. Besides asymptotic convergence results,
we also derive explicit convergence rates for this general
algorithm. It is worth to mention that our convergence  rates depend explicitly on the number of computed projections per iteration.  Moreover, our general framework  generates new algorithms, that are not analyzed in the literature, with possible
better convergence rates than the existing ones. 
%the dependence of the rates on the number of computed projections per iteration $N$, which in the general case
%(no linear regularity) is given by $\mathcal{O}\left(\frac{\gamma_N}{k}\right)$,
%where $\gamma_N = \left(1 - \frac{1}{N}\right)\gamma + \frac{1}{N}$.
%Under the linear regularity assumption, a linear convergence rate $\mathcal{O}\left( \left(1 - \frac{1}{\kappa \gamma_N }\right)^k \right)$
%in the expected distance to the optimal set is obtained,
%where $\kappa$ represents the linear regularity constant.

\noindent $(iii)$ From our convergence analysis it follows that we
can use large step-sizes, besides the usual naturally arisen constant
step-size policy. Thus, we prove theoretically, what is empirically
known in numerical applications for a long time, namely that  these
over-relaxations  accelerate significantly the convergence of
projection methods.

\noindent $(iv)$ Our general random projection algorithm allows to
project simultaneously onto several sets, thus providing great
flexibility in matching the implementation of the algorithms on
the parallel architecture at hand.

\vspace{5pt}

\noindent \textit{Notations}. For given $m \in \mathbb{N}\backslash
\{0\}$, we denote the  set $[m] = \{1, \dots, m \}$. We consider the
space $\rset^n$ composed  by column vectors.  For $x,y \in \rset^n$
denote the scalar product by $\langle x,y \rangle = x^T y$ and the
Euclidean norm by $\|x\|=\sqrt{x^T x}$. We use the notation $x_i$
for the $i$th  component of  the vector $x$ and $e_i$ for the
$i$th column of the identity matrix. The projection operator onto
the closed convex set $X$ is denoted by $\Pi_{X}(\cdot)$ and the distance from a given $x$ to set $X$ is denoted by $\text{dist}_X(x)$. Let $Q \in \rset^{n
\times n}$, then we use notation $Q_i$ for the $i$th row of the
matrix $Q$. The minimal non-zero singular value and the minimal
nonzero eigenvalue of the matrix $Q$ are represented by
$\sigma_{\min}^{\text{nz}}(Q)$ and $\lambda_{\min}^{\text{nz}}(Q)$,
respectively. Similarly, $\sigma_{\max}(Q)$ and $\lambda_{\max}(Q)$
denote the largest singular value and the largest  eigenvalue of the
matrix $Q$, respectively. Also $\|Q\|_F$ denotes its Frobenius norm.

%%%%%%%%%%%%%%%%%%%%%%%%%%%%%%%%%%%%%%%%%%%%%%%%%%%%%%%%%%%%%
%%%%%%%%%%%%%%%%%%%%%%%%%%%%%%%%%%%%%%%%%%%%%%%%%%%%%%%%%%%

\section{Problem formulation}
\noindent In this paper we consider the convex feasibility problem:
\begin{equation}
\label{convexfeas_original_0}
\text{Find} \quad x \in \cX,  %$\eqdef \bigcap_{i=1}^m \cX_i,
\end{equation}
where $\cX \subseteq \rset^n$ is a closed convex set.  We assume
that $\cX$ is nonempty. In general, in most convex feasibility
problems one should seek scalable algorithms with simple iterations
which are able  to find an approximation of a point from the set
$\cX$. For this purpose, we usually assume that $\cX$ can be
represented as the intersections of finitely/infinitely many simple
closed convex sets. Then, a simple and widely  known idea for
solving the convex feasibility problem is to project successively
onto the individual sets in a certain fashion, e.g. cyclic or
random. These projection algorithms are most efficient when  the
projections onto the individual sets are computationally cheap.
However, in many cases, it is difficult to find an explicit
representation of the set $\cX$ as intersection of simple sets. That
is why in the sequel we  consider different relaxations of
\eqref{convexfeas_original_0}, based on several representations for
the individual sets, and we investigate when this relaxations are
exact.

\subsection{Stochastic reformulations}
In many applications the set $\cX$ has  explicit representations,
while in others this set it is not known explicitly. Therefore,
below  we present several representations or approximations for the
set $\cX$. For that, we introduce the concept of  {\em stochastic
approximation} of $\cX$. Given a probability distribution
$\Prob$, we consider a random  variable $S \sim \Prob$ from a
probability space $\Omega$.

\begin{definition}[Stochastic approximation of sets]
\label{def:stoch_approx_sets} For any $S \in \Omega$ let $\cX_S$ be
a random closed convex subset of $\R^n$.  We say that $\cX_S$ is a
{\em stochastic approximation} of $\cX$ if $\cX \subseteq \cX_S$ for
all $S \in \Omega$.
% \textcolor{red}{OR with probability 1, i.e. $\Prob(\cX \subseteq \cX_S) = 1$}.
\end{definition}

\noindent We will henceforth consider stochastic approximation  sets
$\cX_S$ arising as a function of  some random  variable $S$ from a
probability space $(\Omega, \Prob)$.  Therefore, the set $\cX$ may
be represented as an exact  countable/uncountable intersection of
stochastic approximation sets $\cX_S$, that is $\cX = \cap_{S \in \Omega} \cX_S$, or approximated by this intersection, that is $\cX
\subseteq  \cap_{S \in \Omega} \cX_S$. Clearly, having a family of
stochastic approximation sets $(\cX_S)_{S \in \Omega}$, we have the
first relaxation of \eqref{convexfeas_original_0}:
\begin{align}
\label{first_rel} \cX \subseteq \bigcap_{S \in \Omega} \cX_S.
\end{align}
Then, we consider the following convex feasibility problem, which
may be a relaxation of the potentially difficult original problem
\eqref{convexfeas_original_0}:
\begin{equation}
\label{convexfeas_original} \text{Find} \quad x \in  \bigcap_{S \in \Omega} \cX_S.
\end{equation}
In this paper, we propose several  stochastic reformulations of the
convex feasibility problem \eqref{convexfeas_original}.

\begin{enumerate}
\item \textbf{Stochastic fixed point problem.}
\begin{equation}
\label{eq:reform_stoch_fixed_point}\text{Find a fixed point of the
mapping} \quad x \mapsto \mathbf{E}_{S \sim \Prob} \left[ {\Pi_{\cX_{S}}(x)}\right].
\end{equation}

\item \textbf{Stochastic non-smooth optimization problem.}  \begin{equation}
    \label{eq:reform_stoch_opt_indicator}
    \hspace{0.5cm} \text{Minimize} \quad \left\{f(x) \eqdef \mathbf{E}_{S \sim \Prob} \left[ \Ib_{\cX_S}(x) \right] \right\} \quad \text{subject to} \quad {x\in \R^n}.\end{equation}

\item \textbf{Stochastic smooth optimization problem.}  \begin{equation}
    \label{eq:reform_stoch_opt}
 \hspace{2.5cm} \text{Minimize} \quad \left\{F(x)\eqdef  \frac{1}{2} \mathbf{E}_{S \sim \Prob} \left[ \| x - \Pi_{\cX_{S}}(x)\|^2  \right] \right\}
\quad \text{subject to} \quad {x\in \R^n}.\end{equation}

\item \textbf{Stochastic intersection problem.}
\begin{equation}
\label{eq:reform_stoch_intersect}
\hspace{-2cm} \text{Find} \quad x\in \R^n
\quad \text{such that} \quad \Prob(x\in \cX_{S}) = 1.
\end{equation}
\end{enumerate}

\noindent Equivalence of the above reformulations is captured by the
following lemma:
\begin{lemma}[Equivalence]
\label{lem:equivalence} \textcolor{black}{Assume $\cap_{S \sim \Prob}
\cX_S \not = \emptyset$}.  The stochastic reformulations
\eqref{eq:reform_stoch_fixed_point},
\eqref{eq:reform_stoch_opt_indicator},  \eqref{eq:reform_stoch_opt}
and \eqref{eq:reform_stoch_intersect} of the convex feasibility
problem \eqref{convexfeas_original} are equivalent. That is, the set
of fixed points of $x \mapsto \mathbf{E}_{S \sim \Prob} \left[ {\Pi_{\cX_{S}}(x)} \right] $ is equal to
the set of minimizers of the objective functions $f$ or $F$, and
to the set $\{x\;:\; \Prob(x\in \cX_S)=1\}$.  We shall use the
symbol $\cY$ to denote this set.
\end{lemma}

%\begin{proof}  
\noindent \textit{Proof}: An elementary property of the Lebesgue integral  states that if $\phi \geq 0$, then $ \Exp{\phi} = 0 $ if and only if $\phi =0$ almost sure (a.s.). Using this classical result, we can prove the following equivalences:\\

\noindent \eqref{eq:reform_stoch_opt_indicator}$\Leftrightarrow$
\eqref{eq:reform_stoch_intersect}. The  $\Prob$-measurable
function $f_S(x) = \Ib_{\cX_S}(x)$ is non-negative  and thus the set
of minimizers in  \eqref{eq:reform_stoch_opt_indicator} are those
$x$ for which $\mathbf{E}_{S \sim \Prob} \left[ \Ib_{\cX_S}(x) \right] =0$, which is equivalent to
$\Ib_{\cX_S}(x) =0$ a.s., that is $x \in \cX_S$ a.s., or equivalent~to~$\Prob(x\in \cX_{S}) = 1$.\\

\noindent \eqref{eq:reform_stoch_opt}$\Leftrightarrow$
\eqref{eq:reform_stoch_intersect}. The function $F_S(x) = \|
x - \Pi_{\cX_{S}}(x)\|^2$ is non-negative and thus the set of
minimizers in  \eqref{eq:reform_stoch_opt} are those $x$ for which
$\mathbf{E}_{S \sim \Prob} \left[ \| x - \Pi_{\cX_{S}}(x)\|^2 \right] =0$, which is equivalent
to $\| x - \Pi_{\cX_{S}}(x)\| =0$ a.s. or equivalently $ x = \Pi_{\cX_{S}}(x)$ a.s.
or equivalently   $x \in \cX_S$ a.s., or equivalent  to $\Prob(x\in \cX_{S}) = 1$.\\

\noindent
\eqref{eq:reform_stoch_opt}$\Rightarrow$\eqref{eq:reform_stoch_fixed_point}.
Since $\| \mathbf{E}_{S \sim \Prob} \left[ x - \Pi_{\cX_{S}}(x) \right ]\|^2 \leq \mathbf{E}_{S \sim \Prob} \left[ \| x -
\Pi_{\cX_{S}}(x)\|^2 \right ] $, then it follows that the set of minimizers
of \eqref{eq:reform_stoch_opt} are included in the set of fixed
points of  the average projection operator $\Pi(x) =
\mathbf{E}_{S \sim \Prob} \left[ \Pi_{\cX_{S}}(x) \right]$ defined in
\eqref{eq:reform_stoch_fixed_point}. \\

\noindent It remains to prove the other
inclusion
\eqref{eq:reform_stoch_fixed_point}$\Rightarrow$\eqref{eq:reform_stoch_opt}.
Let $x$ be a fixed point of the average projection operator,  that
is $x =\mathbf{E}_{S \sim \Prob} \left[ \Pi_{\cX_{S}}(x) \right] $. Then, for any $z \in \cap_{S \sim \Prob} \cX_S$,  it follows that $z \in \cX_S$ for all $S$ and from
the optimality condition for the projection onto $\cX_S$ we have
$\langle x - \Pi_{\cX_{S}}(x), \Pi_{\cX_{S}}(x)- z   \rangle \geq
0$. This leads to:
\begin{align*}
0 & = \langle  \mathbf{E}_{S \sim \Prob} \left[ x - \Pi_{\cX_{S}}(x) \right] , x- z   \rangle  = \mathbf{E}_{S \sim \Prob} \left[ \langle x - \Pi_{\cX_{S}}(x), x- z   \rangle \right] \\
& =  \mathbf{E}_{S \sim \Prob} \left[ \langle x - \Pi_{\cX_{S}}(x), x - \Pi_{\cX_{S}}(x) + \Pi_{\cX_{S}}(x)- z   \rangle \right] \\
& =  \mathbf{E}_{S \sim \Prob} \left[ \|x - \Pi_{\cX_{S}}(x) \|^2 \right] + \mathbf{E}_{S \sim \Prob} \left[ \underbrace{\langle
x - \Pi_{\cX_{S}}(x), \Pi_{\cX_{S}}(x)- z   \rangle}_{\geq 0} \right],
\end{align*}
for all $z \in \cap_{S \in \Omega} \cX_S$. Thus, sum of two non-negative scalars is zero implies that each term
is zero, that is $ \mathbf{E}_{S \sim \Prob} \left[ \|x - \Pi_{\cX_{S}}(x) \|^2 \right] =0$ and
therefore the set of fixed points of
\eqref{eq:reform_stoch_fixed_point} are included into the set of
minimizers of \eqref{eq:reform_stoch_opt}. \qed
%\end{proof}

%%%%%%%%%%%%%%%%%%%%%%%%%%%%%%%%%%%%%%%%%%%%%%%%%%%%%%%%%%%%%%%%%%%%

\subsection{Discussion}
The proof of Lemma \ref{lem:equivalence} provides several
connections between \eqref{eq:reform_stoch_fixed_point},
\eqref{eq:reform_stoch_opt_indicator},  \eqref{eq:reform_stoch_opt}
and \eqref{eq:reform_stoch_intersect}. There is also an interesting
interpretation between \eqref{eq:reform_stoch_opt_indicator} and
\eqref{eq:reform_stoch_opt}. Notice that for any given  nonempty
closed convex set $Y$, the indicator function $\Ib_{Y}$ is convex,
lower semi-continuous, that is not identically $+ \infty$.
Therefore, the value function:
\[ \frac{1}{2} \norm{x - \Pi_{Y}(x)}^2 = \min\limits_{z \in \rset^n} \Ib_{Y}(z)
 + \frac{1}{2}\norm{z-x}^2 \]
is known  to be well-defined and finite everywhere \cite{BauCom:11}
(Chapter 12). Moreover, the function $ x \mapsto \norm{x -
\Pi_{Y}(x)}^2$ is the Moreau approximation of the non-smooth
indicator function $\Ib_{Y}$, thus it  has Lipschitz continuous
gradient with constant $1$, see \cite{Ned:10}. This implies that the
function $F$ has Lipschitz continuous gradient with constant $L_F =
1$. Observe that the smooth optimization problem
\eqref{eq:reform_stoch_opt} is obtained from the Moreau
approximation $F_S(x) = 1/2 \norm{x - \Pi_{\cX_S}(x)}^2$ of each
indicator function $f_S(x) = \Ib_{\cX_S}(x)$  of  the non-smooth
optimization problem \eqref{eq:reform_stoch_opt_indicator}, that is:
\begin{align}
\label{convexfeas_minimization2} \min_{x \in \rset^n} \; F(x) & \!=\!
\min_{x \in \rset^n} \; \mathbf{E}_{S \sim \Prob} \left[ F_S(x) \right] = \min_{x \in \rset^n} \;
\mathbf{E}_{S \sim \Prob} \left[ \min\limits_{z \in \rset^n} f_{S}(z) \!+\! \frac{1}{2} \norm{z-x}^2 \right]  \\
& \!=\! \min_{x \in \rset^n} \; \mathbf{E}_{S \sim \Prob} \left[ \underbrace{\min\limits_{z \in
\rset^n} \Ib_{\cX_S}(z) + \frac{1}{2}\norm{z-x}^2}_{\norm{x - \Pi_{\cX_S}(x)}^2} \right]. \nonumber
\end{align}
Note that, for general functions $f_S$, there are no connections
between the two problems \eqref{eq:reform_stoch_opt_indicator} and
\eqref{eq:reform_stoch_opt} as expressed in
\eqref{convexfeas_minimization2}. However, for indicator functions
$f_S(x) = \Ib_{\cX_S}(x)$ we have  $\arg \min_x f(x) = \arg \min_x
F(x)$,  according to previous lemma.

\noindent For the convex feasibility problem
\eqref{convexfeas_original}, with $\Omega$ having finite support,
the following basic alternating projection algorithm has been
extensively studied in the literature \cite{Ned:10,GubPol:66}:
\[ (\textbf{B-AP}): \;\;  \text{choose} \; S_k \; \text{cyclic/random \& update} \; x^{k+1} = \Pi_{\cX_{S_k}}(x^k).   \]
The (B-AP) algorithm can be interpreted in several ways depending on the reformulations \eqref{eq:reform_stoch_fixed_point}-\eqref{eq:reform_stoch_intersect}:
\begin{itemize}
\item[1.] For example, when solving the stochastic fixed point problem \eqref{eq:reform_stoch_fixed_point}, we do not have an explicit access to the average projection map $x \to \mathbf{E}_{S \sim \Prob} \left[ \Pi_{\cX_{S}}(x) \right]$.
Instead, we are able to repeatedly sample $S \sim \Prob$  and use
the stochastic projection map $x \to \Pi_{\cX_{S}}(x)$, which leads
to  the random variant of  (B-AP) algorithm.

\item[2.] Since the stochastic  optimization problem
\eqref{eq:reform_stoch_opt_indicator} with $f_S = \Ib_{\cX_S}$:
\[ \min_x f(x) = \mathbf{E}_{S \sim \Prob} \left[ f_S(x) \right ], \]
is non-smooth, then we approximate each indicator function $f_S =
\Ib_{\cX_S}$  with its Moreau approximation $F_S$ and we can apply gradient
method on the resulting expected approximation which leads to the proximal point
method. Since we do not have access to the function $ \mathbf{E}_{S \sim \Prob} \left[ \Ib_{\cX_S}(z) + \frac{1}{2}\norm{z-x}^2 \right] $ for some fixed $x$, but
we can repeatedly sample $S \sim \Prob$ we can apply stochastic
proximal point:
    \[  x^+ = \arg \min_z \Ib_{\cX_S}(z) + \frac{1}{2}\norm{z-x}^2 = \Pi_{\cX_{S}}(x). \]

\item[3.] When solving the stochastic optimization problem  \eqref{eq:reform_stoch_opt}:
\[ \min_x \;  F(x) = \mathbf{E}_{S \sim \Prob} \left[ F_S(x) \right], \]
where
\[ F_S(x) = \frac{1}{2}\norm{x - \Pi_{\cX_S}(x)}^2 \]
we do not have access to the  gradient of $F$:
\[ \nabla F(x) = \mathbf{E}_{S \sim \Prob} \left[ \nabla F_S(x) \right] = \mathbf{E}_{S \sim \Prob} \left[ x - \Pi_{\cX_S}(x) \right]. \]
Instead, we can repeatedly sample $S \sim \Prob$  and receive
unbiased samples of this gradient  at points of interest, that is
$\nabla F_S(x) = x - \Pi_{\cX_S}(x)$. Then, applying the stochastic
gradient method with stepsize $1$ leads to the random variant of (B-AP).

\item[4.] We observe that \eqref{eq:reform_stoch_intersect} can be written equivalently as:
\begin{equation*}
\text{Find} \quad x \in \{x: \;  \Prob(x \in \cX_{S}) = 1 \} := \bigcap_{S \sim \Prob} \cX_S.
\end{equation*}
Then, when solving the previous stochastic intersection problem we
typically do not have explicit access to the stochastic intersection
$\cap_{S \sim \Prob} \cX_S$. Rather, we can sample $S \sim \Prob$
and utilize the simple form of $\cX_S$ to derive (B-AP) algorithm.
\textcolor{black}{Notice that if $\Omega$ is finite/countable, then
the stochastic intersection problem reduces to the standard
intersection problem \eqref{convexfeas_original}.}
\end{itemize}

\noindent However, in Section \ref{sec_spa} we will give a more
general algorithmic framework for solving the four equivalent
problems with larger stepsize and better performances than (B-AP).

\noindent As Lemma~\ref{lem:equivalence} claims, our four stochastic
reformulations are equivalent and they are all   relaxations of the
convex feasibility problem \eqref{convexfeas_original}, that is:
\[ \cap_{S \in \Omega} \cX_S \subseteq \cY. \]
Therefore, for any family
of stochastic approximation sets $(\cX_S)_{S \sim \Prob}$ over a
probability space $(\Omega,\Prob)$, we clearly have:
\begin{equation}
\label{eq:exactness-one-inclusion} \cX \subseteq \bigcap_{S \in \Omega} \cX_S \subseteq \cY.
\end{equation}
Therefore, it is natural to investigate when these inclusions hold
with equality.

%%%%%%%%%%%%%%%%%%%%%%%%%%%%%%%%%%%%%%%%%%%%%%%%%%%%%%%%%%%%%

\section{Exactness}
For simplicity we further redenote $\mathbf{E}_{S \sim \Prob} \left[ \cdot \right]$
with the simpler notation $\mathbf{E} \left[ \cdot \right]$.
From previous discussion we note that for any family of
stochastic approximations $(\cX_S)_{S \sim \Prob}$ over a
probability space $(\Omega,\Prob)$ we trivially have the inclusion
$\cX \subseteq \cY$. If $\cX = \cY$, then the stochastic
reformulations \eqref{eq:reform_stoch_fixed_point},
\eqref{eq:reform_stoch_opt_indicator},  \eqref{eq:reform_stoch_opt}
and \eqref{eq:reform_stoch_intersect}  are equivalent to the convex
feasibility problems \eqref{convexfeas_original_0} and
\eqref{convexfeas_original}. However, this need not be the case, not
without additional assumptions. To see this, consider $\cX =
\bigcap_{i=1}^m \cX_i$, that is finite intersection of closed convex
sets $\cX_i$, and the random set $\cX_S = \cX_1$. Since $\cX
\subseteq \cX_1$, this constitutes a stochastic approximation of
$\cX$, as defined in Definition~\ref{def:stoch_approx_sets}.
However, $\cY=\cX_1$, which is not necessarily equal to $\cX$. In
view of the above, we need to enforce a regularity assumption, which
we call {\em exactness}.

\begin{assumption}[Exactness] Stochastic reformulations
\eqref{eq:reform_stoch_fixed_point},
\eqref{eq:reform_stoch_opt_indicator}, \eqref{eq:reform_stoch_opt}
and \eqref{eq:reform_stoch_intersect}
 of the convex feasibility problems \eqref{convexfeas_original_0} and \eqref{convexfeas_original} are exact. That is, $\cX = \cY$.
\end{assumption}

\noindent In the next result we give a sufficient condition for
exactness:
\begin{lemma} The following statement hold:
If there exists $\kappa < \infty$ such that  the following inequality
(a.k.a. ``linear regularity property'') holds for all $x \in \R^n$:
\begin{equation}
\label{linreg}
\dist_{\cX}^2(x) \le \kappa \; \Exp{ \dist_{\cX_S}^2(x)},
\end{equation}
then $\cX = \cY$ (i.e., exactness holds).
\end{lemma}

%\begin{proof}
\noindent \textit{Proof}:  The set $\cY$ of optimal points  of the stochastic smooth optimization problem \eqref{eq:reform_stoch_opt} satisfies:  $F(x)
=0$ for all $x \in \cY$. Moreover, the relation  $F(x) = \Exp{ \| x -
\Pi_{\cX_{S}}(x)\|^2  } = \Exp{\dist_{\cX_S}^2(x)}$ holds. Therefore, for
any  $x\in \cY$ we have $\Exp{\dist_{\cX_S}^2(x)} = 0$. From
\eqref{linreg} we conclude that $\dist_{\cX}^2(x)=0$, which means
that $x\in \cX$. Combined  with \eqref{eq:exactness-one-inclusion},
this implies  that $\cX = \cY$ holds. \qed
%\end{proof}

\noindent Since  $\dist_{\cX_S}(x) \leq \dist_{\cX}(x)$ it follows
immediately from \eqref{linreg} that $\kappa \geq 1$. The
feasibility problem is ill-conditioned when $\kappa$ is large.
\textcolor{black}{Notice that linear regularity is a very conservative
condition for exactness. We can see that  if $\Omega$ is
finite/countable, then the stochastic intersection problem reduces
to the standard intersection problem \eqref{convexfeas_original},
i.e. we have exactness.}  Note that linear regularity property does
not hold for any collection of closed convex sets as the following
example  shows:

\begin{example}  Let $\cX_1 =\{x: |x_1|^p \leq x_2\}$ with $p>1$, and
$\cX_2=\{x: x_2=0\}$. These two sets are convex and $\cX = \cX_1
\cap \cX_2 = \{0\}$.  Then, for any $x \in \cX_1$, satisfying
$|x_1|^p = x_2$, we have:
\[ \dist_{\cX}^2(x) = x_1^2 + x_2^2 \quad \text{and} \quad  \dist_{\cX_1}^2(x) + \dist_{\cX_2}^2(x) = x_2^2. \]
Then, clearly there is no finite $\kappa >0$ such that:
\[   x_1^2 + x_2^2 \leq \kappa x_2^2 \quad \forall |x_1|^p = x_2, \; x_1\geq 0, \]
since by replacing $x_2$ and obtaining
\[  x_1^2 + x_1^{2p} \leq \kappa x_1^{2p} \Rightarrow \frac{1}{x_1^{2p-2}} +1 \leq \kappa, \]
we can take $x_1$ very small (close to zero). \qed
\end{example}

\noindent Notice that linear regularity is related to Slater's condition, as discussed in \cite{BecTeb:03}.  Moreover, this property  is directly
related to the stochastic formulations \eqref{eq:reform_stoch_fixed_point}-\eqref{eq:reform_stoch_intersect},
as we will show below.

%%%%%%%%%%%%%%%%%%%%%%%%%%%%%%%%%%%%%%%%%%%%%%%%%%%%%%

\subsection{Properties of the smooth function $F$}
\label{sec_F}
If we consider    the smooth stochastic optimization
problem \eqref{eq:reform_stoch_opt}, then we have the following
important relation:
\begin{align}
\label{relation_linSGD} F(x) = \frac{1}{2} \Exp{\| \nabla F_S(x)\|^2
} \quad \forall x \in \rset^n,
\end{align}
since we recall that $\nabla F_S(x) = x - \Pi_{\cX_S}(x)$. Moreover,
the linear regularity property \eqref{linreg} is equivalent with the
quadratic functional growth condition on $F$ introduced in
\cite{NecNes:15}, which was defined as a relaxation of strong
convexity. Indeed, under the exactness assumption, we have $\cX= \cY
= \arg \min_x F(x)$ and the optimal value $F^*=0$. Moreover, we have
$ F(x) = \Exp {\dist_{\cX_S}^2(x) }$.  Then, the property
\eqref{linreg} can be rewritten equivalently as:
\begin{align}
\label{qg} F(x) -F^*   \geq \frac{1}{2 \kappa} \| x - \Pi_{\cX}(x)
\|^2 \quad \forall x \in \rset^n,
\end{align}
which is exactly the definition of the quadratic functional growth
condition introduced in \cite{NecNes:15}. Typically, the standard
assumption for proving linear  convergence of first order methods
for smooth convex optimization is the strong convexity of the
objective function, an assumption which does not hold for many
practical applications, including the one presented in this paper.
In \cite{NecNes:15} it has been proved that we can still achieve
linear convergence rates of several first order methods for solving
smooth non-strongly convex constrained optimization problems, i.e.
involving an objective function with a Lipschitz continuous gradient
that satisfies the relaxed strong convexity condition \eqref{qg}.
Moreover, in \cite{NecNes:15} it has been shown that the quadratic
functional growth condition \eqref{qg} is equivalent with the
so-called error bound condition for unconstrained problem
\eqref{eq:reform_stoch_opt}.

\noindent Further, let $ \Uconst \geq 0$ be the smallest constant
satisfying the inequality:
\begin{equation}\label{eq:gamma0}
\left\|\Exp{x-\Pi_{\cX_{S}}(x)}\right\|^2 \leq \Uconst \cdot \Exp{
\left\| x - \Pi_{\cX_{S}}(x) \right\|^2 } \qquad \forall x\in \R^n.
\end{equation}
By Jensen's inequality,  $\Uconst \leq 1$. However, for specific
sets and distributions $\Prob$, it is possible for $\Uconst$ to be
strictly smaller than $1$, as the following examples show. For
example, we can consider finding a solution of a linear system
$\cX=\{x: A x = b \}$, where $A \in \rset^{m \times n}$. For this
set we can easily construct stochastic approximations sets $\cX_S =
\{x: \; S^T A x = S^Tb \}$ taking any matrix  $S \in \rset^{m \times
q}$. Clearly, for any matrix $S$ we have $\cX \subseteq \cX_S$.
Then, we have the following characterization for $\Uconst$:

\begin{theorem}
\label{th1_Ules1} Let us consider finding a solution of the linear
system $\cX=\{x: A x = b \}$, where $A \in \rset^{m \times n}$.
Further, let us consider the stochastic approximation sets $\cX_S =
\{x: \; S^T A x = S^Tb \}$, where $S \in \Omega = \rset^{m \times
q}$ and a probability distribution $\Prob$ on $\Omega$. Then,
\eqref{eq:gamma0} holds with:
\[  \Uconst = \lambda_\text{max}(A^T \Exp{S(S^TA A^T S)^{\dagger}S^T} A) \leq 1.  \]
\end{theorem}

%\begin{proof}
\noindent \textit{Proof}: Clearly, for $x$ satisfying $Ax=b$ the inequality \eqref{eq:gamma0}
holds for any $\Uconst \leq 1$. It remains to prove for $x$
satisfying $Ax-b \not = 0$. However,  since $\cX_S = \{x: \; S^T A x
= S^Tb \}$, then the projection of $x$ onto $\cX_S$ can be computed
explicitly:
\[  \Pi_{\cX_{S}}(x) =x - A^T S (S^TAA^TS)^{\dagger}S^T(A x - b)  \]
and the relation we need to prove becomes as follows:
\[  \left\| \Exp{A^TS (S^TAA^TS)^{\dagger}S^T(A x - b)} \right\|^2 \leq \Uconst \Exp{ \| A^TS (S^TAA^TS)^{\dagger} S^T(A x - b) \|^2 }.    \]  Using
the standard properties of the pseudoinverse, that is $Q^\dagger Q
Q^\dagger = Q^\dagger$ for any matrix $Q$, the previous relation is
equivalent to:
\[  \left\| A^T \Exp{S (S^TAA^TS)^{\dagger}S^T}(A x - b) \right\|^2 \leq
\Uconst (Ax-b)^T \Exp{S (S^TAA^TS)^{\dagger}S^T}(Ax-b). \] For
simplicity, let us   denote  $E = \Exp{S (S^TAA^TS)^{-1}S^T}$. Then
$E$ is a positive semidefinite matrix and thus there exists
$E^{1/2}$.  Clearly, for $Ax-b \in \text{Null}(E)$ the previous
inequality holds for any $\Uconst$. Therefore, $\Uconst$ is defined as:
\begin{align*}
\Uconst & = \max_{x: Ax -b \not \in \text{Null}(E)}
\frac{\|A^T E (Ax-b)\|^2}{(Ax-b)^T E (Ax-b)}  \\
& = \max_{x: Ax -b \not \in \text{Null}(E^{1/2})} \frac{\|A^T
E^{1/2} E^{1/2} (Ax-b) \|^2}{\| E^{1/2} (Ax-b)\|^2} \\
& = \max_{z \not =0} \frac{\|A^T E^{1/2} z \|^2}{\| z\|^2}
 = \sigma_\text{max}^2(A^T E^{1/2}) = \lambda_\text{max}(A^T E A).
\end{align*}
Therefore, we have $\Uconst = \lambda_\text{max}(A^T \Exp{S
(S^TAA^TS)^{\dagger}S^T} A)$. Since the function $W \mapsto
\lambda_\text{max}(W)$ is convex over the space of positive
semidefinite matrices, then using Jensen's inequality we have:
\begin{align*}
\lambda_\text{max}(A^T \Exp{S(S^TAA^TS)^{\dagger}S^T} A) \leq \Exp{
\lambda_\text{max}(A^T S (S^TAA^TS)^{\dagger}S^T A)}.
\end{align*}
Furthermore, the matrix $P_S = A^T S (S^TAA^TS)^{\dagger}S^T A$ is
idempotent, that is $P_S^2 =P_S$. Therefore, all the eigenvalues of
$P_S$ are either $0$ or $1$. Then, we get:
\[ \Uconst = \lambda_\text{max}(A^T \Exp{S
(S^TAA^TS)^{\dagger}S^T} A) \leq  \Exp{ \lambda_\text{max}(A^T S
(S^TAA^TS)^{\dagger}S^T A)} \leq 1, \] which proves the statement of
the theorem. \qed
%\end{proof}

\noindent Based on the previous theorem we can prove that  for
particular choices of the probability distribution $\Prob$ we have
$\Uconst < 1$, see e.g. the next corollary:

\begin{corollary}
\label{cor1_Ules1} Let us consider finding a solution of the linear
system $\cX=\{x: A x = b \}$, where $A \in \rset^{m \times n}$
having  $\text{rang}(A) \geq 2$.  Further, let us consider
$\Omega=\{e_1, \cdots, e_m\}$, the standard basis of $\rset^m$, and
the corresponding stochastic approximation sets $ \cX_{e_i} = \{x:
\; A_i^T x = b_i \}$ for all $i \in [m]$.  Then, for two choices of
the probability distribution $\Prob$ on $\Omega$, inequality
\eqref{eq:gamma0} holds with:
\begin{align}
\label{cor1_Ules2}
\Uconst =
\begin{cases}
\frac{\lambda_\text{max} \left(A^T A \right)}{\|A\|_F^2} & \quad
\text{if} \quad \Prob(S=e_i) = \frac{\norm{A_i}^2}{\norm{A}^2_F}\\
\frac{\lambda_\text{max} \left(A^T D A \right)}{m} & \quad \text{if}
\quad \Prob(S=e_i) = \frac{1}{m}
\end{cases}
\quad <1,
\end{align}
where the diagonal matrix $D \eqdef \text{diag}(\|A_1\|^{-2},
\cdots, \|A_m\|^{-2})$.
\end{corollary}

%\begin{proof}
\noindent \textit{Proof}: In this case we have the following expression: $$S
(S^TAA^TS)^{\dagger}S^T = e_i (e_i^TA A^T e_i)^{\dagger} e_i^T =
\frac{1}{\|A_i\|^{2}} e_i e_i^T. $$ Then, from Theorem
\ref{th1_Ules1} we get for  probability distribution
$\Prob(S=e_i)~=~\frac{\norm{A_i}^2}{\norm{A}^2_F}$:
\begin{align*}
\Uconst & = \lambda_\text{max}\left(A^T \Exp{ \frac{1}{\|A_i\|^{2}} e_i e_i^T}
A \right)  = \lambda_\text{max}\left(A^T \sum_{i=1}^m  \frac{\|A_i\|^2}{\|A\|^2_F} \frac{1}{\|A_i\|^{2}} e_i e_i^T A \right) \\
&= \lambda_\text{max} \left( \frac{A^T A}{\|A\|_F^2} \right) =  \lambda_\text{max} \left( \frac{A A^T}{\|A\|_F^2} \right).
\end{align*}
In the last equality we used the fact that the maximum eigenvalues
of the matrices $A^T A$ and $A A^T$ coincides. But we can easily see
that the trace of the matrix $\frac{A A^T}{\|A\|_F^2}$ is equal to
$1$ and thus:
\[ \sum_{i=1}^m \lambda_i \left( \frac{A A^T}{\|A\|_F^2} \right) = \text{Trace}\left( \frac{A A^T}{\|A\|_F^2} \right) =1. \]  Therefore, if
$\text{rang}(A) \geq 2$, then $\Uconst = \lambda_\text{max} \left(
\frac{A A^T}{\|A\|_F^2} \right) < 1$. Similarly,  from Theorem
\ref{th1_Ules1} we obtain for the uniform  probability distribution
$\Prob(S=e_i) = \frac{1}{m}$:
\begin{align*}
\Uconst & = \lambda_\text{max}\left(A^T \Exp{ \frac{1}{\|A_i\|^{2}}
e_i e_i^T} A \right)  = \lambda_\text{max}\left(A^T \sum_{i=1}^m  \frac{1}{m} \frac{1}{\|A_i\|^{2}}
e_i e_i^T A \right) \\
&= \lambda_\text{max} \left( \frac{A^T D A}{m} \right) =
\lambda_\text{max} \left( \frac{A A^T D}{m} \right),
\end{align*}
where $D = \text{diag}(\|A_1\|^{-2}, \cdots, \|A_m\|^{-2})$ and we
used  the fact that the sets of nonzero eigenvalues of the matrices
$U V$ and $VU$ are the same for any two matrices $U$ and $V$ of
appropriate dimensions, in particular $U = A^T D$ and $V = A$.
Moreover, the trace of the matrix $\frac{A A^T D}{m}$ is equal to
$1$ and thus:
\[ \sum_{i=1}^m \lambda_i \left( \frac{A A^T D}{m} \right) = \text{Trace}\left( \frac{A A^T D}{m} \right)
=1. \]  If $\text{rang}(A) \geq 2$, then $\Uconst =
\lambda_\text{max} \left( \frac{A A^T D}{m} \right) < 1$ for uniform
distribution. \qed
% \end{proof}

\noindent For systems of linear inequalities we can obtain  similar
statements. For example, we can consider finding a feasible point
for a system of  linear inequalities $\cX=\{x: A x \leq b \}$, where
$A \in \rset^{m \times n}$. For this set we can easily construct
stochastic approximations sets $\cX_S = \{x: \; S^T A x \leq S^Tb
\}$, where $S$ is a vector with nonnegative entries, i.e.   $S \in
\rset^{m}_+$. Clearly, if the vector $S$ has nonnegative entries, we
have $\cX \subseteq \cX_S$. Then, we have the following
characterization for $\Uconst$:

\begin{theorem}
\label{th2_Ules1} Let us consider finding a solution of a system of
linear inequalities $\cX=\{x: A x \leq b \}$, where $A \in \rset^{m
\times n}$. Further, let us consider the stochastic approximation
sets $\cX_S = \{x: \; S^T A x \leq S^Tb \}$, where $S \in \Omega =
\rset^{m}_+$ and a probability distribution $\Prob$ on $\Omega$.
Then, \eqref{eq:gamma0} holds with:
\[  \Uconst = \lambda_\text{max} \left( A^T \Exp{S(S^TA A^T S)^{-1}S} A \right) \leq 1.  \]
\end{theorem}

%\begin{proof}
\noindent \textit{Proof}: Clearly, for $x$ satisfying $Ax \leq b$ the inequality
\eqref{eq:gamma0} holds for any $\Uconst \leq 1$. It remains to
prove for $x$ satisfying $Ax \not \leq b$. However,  since $\cX_S =
\{x: \; S^T A x \leq  S^Tb \}$, then the projection of $x$ onto
$\cX_S$ can be computed explicitly:
\[  \Pi_{\cX_{S}}(x) =x -  \frac{\max (0,S^T(A x - b))}{\|A^TS\|^2} A^TS =  x -  \frac{\Pi_+ (S^T(A x - b))}{\|A^TS\|^2} A^TS \]
and the relation we need to prove becomes as follows:
\begin{align*}
& \left\| \Exp{A^TS (S^TAA^TS)^{-1} \Pi_+ (S^T(A x - b))} \right\|^2 \leq \Uconst \Exp{ \| A^TS (S^TAA^TS)^{-1} \Pi_+ (S^T(A x - b)) \|^2}
\end{align*}
or equivalently
\begin{align*}
& \left\| A^T \Exp{S (S^TAA^TS)^{-1} \Pi_+ (S^T(A x - b))} \right\|^2 \leq \Uconst \Exp{ \Pi_+(S^T(A x - b)) (S^TAA^TS)^{-1} \Pi_+(S^T(A x - b))}.
\end{align*}
Moreover,  if  we define the event ${\cal I}(x) =\{ S \in
\Omega: S^T(A x - b) > 0 \}$, then the previous relation can be written as follows:
\begin{align*}
& \left\| A^T \left( \int_{{\cal I}(x)}  S (S^TAA^TS)^{-1} S^T dP \right) (A x - b) \right\|^2  \leq \Uconst (A x - b)^T  \left( \int_{{\cal I}(x)}  S (S^TAA^TS)^{-1} S^T dP \right) (A x - b).
\end{align*}
Let us define  $E(x) = \int_{{\cal I}(x)}  S (S^TAA^TS)^{-1} S dP$ and $E = \int_{\Omega}  S (S^TAA^TS)^{-1} S dP$. Then both matrices  are positive semidefinite and $E(x) \preceq E$ for all $x$ such that $Ax \not \leq b$. It follows that $\Uconst$ is an upper bound on the following function:
\begin{align*}
{\cal R}(x) =  \frac{\|A^T E(x) (Ax-b)\|^2}{(Ax-b)^T E(x) (Ax-b)} \leq \Uconst \quad \forall x:\; Ax \not \leq b.
\end{align*}
However, it is easy to find an upper bound for this function ${\cal R}(x)$ for each fixed $x$, namely:
\[  {\cal R}(x) \leq  \lambda_\text{max}(A^T E(x) A) \quad \forall x:\; Ax \not \leq b. \]
Since $E(x) \preceq E$, then  $A^T E(x) A \preceq A^T E A$ and
consequently $\lambda_\text{max}(A^T E(x) A) \leq
\lambda_\text{max}(A^T E A)$. Moreover, there exists  $x$ such that
${\cal I}(x) = \Omega$. Thus,  we have:
\[ \Uconst = \lambda_\text{max}(A^T E A) =  \lambda_\text{max} \left(A^T \Exp{S
(S^TAA^TS)^{-1}S^T} A \right). \] Since the function $W \mapsto
\lambda_\text{max}(W)$ is convex over the space of positive
semidefinite matrices, then using Jensen's inequality we have:
\begin{align*}
\lambda_\text{max} \left(A^T \Exp{S(S^TAA^TS)^{-1}S^T} A \right) \leq \Exp{
\lambda_\text{max} \left(A^T S (S^TAA^TS)^{-1}S^T A \right)}.
\end{align*}
Furthermore, the matrix $P_S = A^T S (S^TAA^TS)^{-1}S^T A$ is
idempotent, that is $P_S^2 =P_S$. Therefore, all the eigenvalues of
$P_S$ are either $0$ or $1$. Then, we get:
\[ \Uconst \!=\! \lambda_\text{max}  \!\left(A^T \Exp{S
(S^TAA^TS)^{-1}S^T} A \right) \!\leq\!  \Exp{ \lambda_\text{max} \! \left(A^T S
(S^TAA^TS)^{-1}S^T A \right)} \!\leq\! 1, \] which proves the statement of
the theorem. \qed
%\end{proof}

\noindent Based on the previous theorem we can prove that  for
particular choices of the probability distribution $\Prob$ we have
$\Uconst < 1$, see e.g. the next corollary:

\begin{corollary}
\label{cor2_Ules1} Let us consider solving a  system of linear
inequalities  $\cX=\{x: A x  \leq b \}$, where $A \in \rset^{m
\times n}$ having  $\text{rang}(A) \geq 2$.  Further, let us
consider $\Omega=\{e_1, \cdots, e_m\}$, the standard basis of
$\rset^m$, and the corresponding stochastic approximation sets $
\cX_{e_i} = \{x: \; A_i^T x \leq  b_i \}$ for all $i \in [m]$. Then,
for two choices of the probability distribution $\Prob$ on $\Omega$,
inequality \eqref{eq:gamma0} holds with:
\begin{align}
\label{cor2_Ules2} \Uconst =
\begin{cases}
\frac{\lambda_\text{max} \left(A^T A \right)}{\|A\|_F^2} & \quad
\text{if} \quad \Prob(S=e_i) = \frac{\norm{A_i}^2}{\norm{A}^2_F}\\
\frac{\lambda_\text{max} \left(A^T D A \right)}{m} & \quad \text{if}
\quad \Prob(S=e_i) = \frac{1}{m}
\end{cases}
\quad <1,
\end{align}
where the diagonal matrix $D \eqdef \text{diag}(\|A_1\|^{-2},
\cdots, \|A_m\|^{-2})$.
\end{corollary}

%\begin{proof} 
\noindent \textit{Proof}: The proof is similar to the one given in Corollary
\ref{cor1_Ules1}. \qed
%\end{proof}

\noindent The reader can easily find   other examples of convex
feasibility problems with $\Uconst <1 $. The linear regularity
inequality \eqref{linreg} and the Jensen type inequality
\eqref{eq:gamma0} impose strong conditions on the shape of the
function $F$:
\begin{theorem}
\label{th_shapeF} Let the linear regularity condition \eqref{linreg}
hold. Then, the following bounds are valid for the smooth objective
function $F$:
\begin{align}
\label{qg_Lip} \frac{1}{2 \kappa} \| x - \Pi_{\cX}(x) \|^2 \leq F(x)
-F^* \leq \frac{\Uconst}{2} \| x - \Pi_{\cX}(x) \|^2 \quad \forall x
\in \rset^n,
\end{align}
and their dual formulations
\begin{align}
\label{qg_Lip2} \frac{1}{2 \Uconst} \|\nabla F(x)\|^2  \leq
 F(x) - F^* \leq \frac{\kappa}{2} \|\nabla F(x)\|^2
\quad \forall x \in \rset^n.
\end{align}
\end{theorem}

%\begin{proof}
\noindent \textit{Proof}: Under the linear regularity condition \eqref{linreg} we have \eqref{qg}, which represents the left hand side inequality in
\eqref{qg_Lip}.  For proving the right hand side inequality in
\eqref{qg_Lip} we use a well-known property of the projection:
\begin{align}
\label{proj} \| x - \Pi_{\cX_S} (x) \|^2 \leq \|x - z\|^2 - \|
\Pi_{\cX_S} (x) - z \|^2 \quad \forall z \in \cX_S.
\end{align}
Then, using that $\Pi_{\cX}(x) \in \cX_S$ we have:
\begin{align*}
\Exp { \| x - \Pi_{\cX_S}(x) \|^2 }  & = \| x - \Pi_{\cX}(x) \|^2 + \Exp { \| \Pi_{\cX}(x) -
\Pi_{\cX_S}(x) \|^2 } + 2  \langle x -
\Pi_{\cX}(x), \Exp { \Pi_{\cX}(x) - \Pi_{\cX_S}(x) } \rangle  \\
& \overset{\eqref{proj}}{\leq} 2 \| x - \Pi_{\cX}(x) \|^2 - \Exp { \| x - \Pi_{\cX_S}(x)
\|^2 } + 2  \langle x - \Pi_{\cX}(x), \Exp { \Pi_{\cX}(x) -
\Pi_{\cX_S}(x) } \rangle \\
& = - \Exp { \| x - \Pi_{\cX_S}(x) \|^2 } + 2  \langle x -
\Pi_{\cX}(x), \Exp { x - \Pi_{\cX_S}(x) } \rangle,
\end{align*}
where in the first inequality we used $\eqref{proj}$.
In conclusion, we get:
\begin{align}
\label{proj_imp}  \Exp { \| x - \Pi_{\cX_S}(x) \|^2 }  & \leq
\langle x - \Pi_{\cX}(x), \Exp { x - \Pi_{\cX_S}(x) } \rangle.
\end{align}
Furthermore, using  Cauchy-Schwartz  inequality and
\eqref{eq:gamma0} in \eqref{proj_imp} we get:
\begin{align*}
\Exp { \| x - \Pi_{\cX_S}(x) \|^2 }  & \leq   \| x - \Pi_{\cX}(x) \|
\| \Exp{ x - \Pi_{\cX_S}(x) } \| \\
& \leq \| x - \Pi_{\cX}(x) \| \sqrt{ \Uconst \Exp { \| x -
\Pi_{\cX_S}(x) \|^2 } }
\end{align*}
which leads to: \[ F(x) - F^* \leq \frac{\Uconst}{2} \| x -
\Pi_{\cX}(x) \|^2,   \] that is, the right hand side inequality in
\eqref{qg_Lip}  holds. This proves the first statement of the
theorem, i.e. \eqref{qg_Lip}.

\noindent For proving the second statement, \eqref{qg_Lip2}, we
first notice that  since the Jensen type inequality
\eqref{eq:gamma0} always holds for some $\Uconst \leq 1$ and using
the expression of $F$ and that $F^*=0$, then we can easily find the
left hand side inequality in~\eqref{qg_Lip2}:
\begin{align}
\label{lip} \frac{1}{2} \| \nabla F(x) \|^2 \leq \Uconst (F(x) -
F^*) \quad \forall x \in \rset^n.
\end{align}
Then, combining  \eqref{linreg} and \eqref{proj_imp} and using
Cauchy-Schwartz inequality,  we get:
\begin{align*}
 \Exp { \| x - \Pi_{\cX_S}(x) \|^2 } & \leq    \| x - \Pi_{\cX}(x)
\|  \| \Exp { x - \Pi_{\cX_S}(x) } \| \\
& \leq  \sqrt{ \kappa  \Exp { \| x - \Pi_{\cX_S}(x) \|^2 } } \| \Exp
{ x - \Pi_{\cX_S}(x) } \|,
\end{align*}
which leads to
\begin{align*}
\frac{1}{\kappa} (F(x) - F^*) \leq  \frac{1}{2} \|\nabla F(x)\|^2.
\end{align*}
Combining the previous inequality with \eqref{lip} we obtain the
second statement of the theorem, i.e. \eqref{qg_Lip2}.  \qed
%\end{proof}

\noindent Theorem \ref{th_shapeF} states that $F$ is strongly convex
with constant $\frac{1}{k}$ and has Lipschitz continuous gradient
with constant $\Uconst$ when restricted along any segment $[x,
\Pi_{\cX}(x)]$. Indeed, since $\nabla F(\Pi_{\cX}(x)) =0$, then from
\eqref{qg_Lip}-\eqref{qg_Lip2} we obtain:
\begin{align*}
\frac{1}{2 \kappa} \| x - \Pi_{\cX}(x) \|^2 & + \langle \nabla
F(\Pi_{\cX}(x)), x - \Pi_{\cX}(x) \rangle + F^* \\
& \leq F(x) \leq
\frac{\Uconst}{2} \| x - \Pi_{\cX}(x) \|^2 + \langle \nabla
F(\Pi_{\cX}(x)), x - \Pi_{\cX}(x) \rangle + F^*
\end{align*}
which are exactly the strong convexity condition and the Lipschitz
continuity condition, respectively,  along any segment $[x,
\Pi_{\cX}(x)]$, see \cite{NecNes:15} for more details. It follows
that $\kappa \Uconst \geq 1$ and $\kappa \Uconst$ represents  the
condition number of the convex feasibility problem
\eqref{convexfeas_original}. Note that $F$ has global  Lipschitz
continuous gradient with constant $L_F=1$.

%%%%%%%%%%%%%%%%%%%%%%%%%%%%%%%%%%%%%%%%%%%%%%%%%%%%%%%%%%%%%%

\subsection{Properties of the operator $\Pi = \Exp{\Pi_{\cX_S}}$}
It is well-known that the projection operator is firmly
nonexpansive:
\[  \langle \Pi_{\cX_S}(x) - \Pi_{\cX_S}(y), x - y \rangle \geq \| \Pi_{\cX_S}(x) -
\Pi_{\cX_S}(y) \|^2 \quad \forall x,y \in \rset^n. \] Taking the
expectation in the previous relation, we get that average projection
operator $\Pi(x) = \Exp{\Pi_{\cX_S}(x)}$ is also firmly
nonexpansive:
\begin{align} \label{oper_firm}
\langle \Exp{\Pi_{\cX_S}(x)} -  \Exp{\Pi_{\cX_S}(y)}, x - y \rangle
\geq \Exp { \| \Pi_{\cX_S}(x) - \Pi_{\cX_S}(y) \|^2 } \geq  \| \Exp{\Pi_{\cX_S}(x)} - \Exp{\Pi_{\cX_S}(y)} \|^2
\end{align}
for all $x,y \in \rset^n$. Similar to Theorem \ref{th_shapeF} we can
derive some contraction inequalities for the average projection
operator $\Pi$.
\begin{theorem}
\label{th_shapeT} Let the linear regularity condition \eqref{linreg}
hold. Then, the following bounds are valid for the average
projection operator $\Pi(x) = \Exp{\Pi_{\cX_S}(x)}$:
\begin{align} \label{in_oper}
\textcolor{black}{\left( 1  - \Uconst \right)} \|x - x^* \|^2
\leq \langle \Pi(x) - \Pi(x^*), x - x^* \rangle \leq   \left( 1  -
\frac{1}{\kappa} \right) \|x - x^* \|^2
\end{align}
for all $x \in \rset^n$ and the corresponding  fixed point  $x^* = \Pi_{\cX}(x)$.
\end{theorem}

%\begin{proof}
\noindent \textit{Proof}: In order to prove the right hand side inequality, we choose in \eqref{oper_firm} the fixed point $y = \Pi_{\cX}(x)$,  which leads to:
\begin{align*}
&  \langle \Exp{ \Pi_{\cX_S}(x)} - \Pi_{\cX}(x), x - \Pi_{\cX}(x) \rangle  \geq \Exp { \| \Pi_{\cX_S}(x) - \Pi_{\cX}(x) \|^2 } \\
& = \Exp{ \| \Pi_{\cX_S}(x) - x \|^2 } - \|x - \Pi_{\cX}(x) \|^2 + 2
\langle \Exp{ \Pi_{\cX_S}(x)} - \Pi_{\cX}(x), x -  \Pi_{\cX}(x)
\rangle,
\end{align*}
which combined with \eqref{linreg} leads to
\begin{align*}
\langle \Exp{ \Pi_{\cX_S}(x)} - \Pi_{\cX}(x), x -  \Pi_{\cX}(x)
\rangle & \leq \|x - \Pi_{\cX}(x) \|^2 - \Exp{ \| \Pi_{\cX_S}(x) - x
\|^2 } \\
& \overset{\eqref{linreg}}{\leq}  \left( 1  - \frac{1}{\kappa}
\right) \|x - \Pi_{\cX}(x) \|^2.
\end{align*}
For the left hand side inequality we proceed as follows:
\begin{align*}
 \langle \Exp{ \Pi_{\cX_S}(x)} - \Pi_{\cX}(x), x -  \Pi_{\cX}(x) \rangle
 & =  \|x - \Pi_{\cX}(x) \|^2  + \langle \Exp{ \Pi_{\cX_S}(x)} - x, x
- \Pi_{\cX}(x) \rangle \\
& \geq \|x - \Pi_{\cX}(x) \|^2  -  \| x - \Pi_{\cX}(x)\|  \|  \Exp{
\Pi_{\cX_S}(x)} - x\| \\
& \geq  \|x - \Pi_{\cX}(x) \|^2  - \| x - \Pi_{\cX}(x)\| \sqrt{
\Uconst \Exp{ \| \Pi_{\cX_S}(x) - x\|^2 }} \\
&  =  \|x - \Pi_{\cX}(x) \|^2  - \| x - \Pi_{\cX}(x)\| \sqrt{
2 \Uconst (F(x) - F^*)} \\
%&  \overset{\eqref{qg_Lip}}{\ge}  \|x - \Pi_{\cX}(x) \|^2  - \| x - \Pi_{\cX}(x)\|  \\
&  \overset{\eqref{qg_Lip}}{\ge}  \|x - \Pi_{\cX}(x) \|^2  -  \| x - \Pi_{\cX}(x)\|
\Uconst \sqrt{\|x - \Pi_{\cX}(x) \|^2} \\
& = (1 - \Uconst) \|x - \Pi_{\cX}(x) \|^2,
\end{align*}
where in the first inequality we used Cauchy-Schwartz inequality, in
the second inequality we used \eqref{eq:gamma0} and in the third
inequality we used relation \eqref{qg_Lip}.  \qed
% \end{proof}

\noindent Theorem \ref{th_shapeT}  shows that the operator $\Pi$ is a
contraction with contraction constant $c = 1 - \frac{1}{k} <1$ when
restricted along any segment $[x, \Pi_{\cX}(x)]$.

%%%%%%%%%%%%%%%%%%%%%%%%%%%%%%%%%%%%%%%%%%%%%%%%%%%%%%%%%%%%

\section{Examples: finite intersection}
\label{sec_finite_inters}
We consider $X$ represented
as the intersection of a finite family of convex sets:
\[  \cX=\bigcap_{i=1}^m \cX_i, \]
where $\cX_i$ are nonempty closed  convex sets. We also assume that
$\cX \not = \emptyset$.  In several papers,  such as
\cite{Ned:10, BauBor:96}, the authors introduced a \textit{linear
regularity property} for the set $\cX = \bigcap_{i=1}^m \cX_i$.  That is,
there exists $\kappa_{\max} < \infty $ such that:
\begin{equation}
\label{linreg_finite}
\dist_{\cX}^2(x) \le \kappa_{\max} \; \max\limits_{i \in [m]}
\dist_{\cX_i}^2(x) \qquad  \forall x \in \rset^n.
\end{equation}
Based on this condition, linear convergence rate, depending on the
constant $\kappa_{\max}$, has been derived for  the alternating
projection algorithm (B-AP). Note that our definition of linear
regularity \eqref{linreg}  extends the one given in
\eqref{linreg_finite}  for finite intersection to the more general
convex feasibility problem \eqref{convexfeas_original}. More
precisely, in order to show linear convergence for our general
algorithmic framework introduced in this paper, we require the
linear regularity property for the set $\cX = \cap_{S \in \Omega}
\cX_S$ defined in \eqref{linreg}.  For a uniform probability over
the set $\Omega= [m] \eqdef \{1,2,\dots,m\}$ we have:
\begin{align*}
\dist_{\cX}^2(x) & \le \kappa_{\max} \; \max\limits_{i \in [m]}
\dist_{\cX_i}^2(x)  \leq  \kappa_{\max} \; \sum_{i=1}^m \dist_{\cX_i}^2(x) \\
& = m \cdot \kappa_{\max} \Exp{ \dist_{\cX_S}^2(x)}.
\end{align*}
This shows that:
\[ \kappa  \leq  m \cdot \kappa_{\max}.  \]
Thus, condition \eqref{linreg_finite} is a relaxation of our more
general condition \eqref{linreg},  also analyzed in \cite{Ned:11}.
Further, we analyze this property \eqref{linreg} and estimate the
constant $\kappa$ for several representative  cases of stochastic
approximation sets for~$\cX$.

%%%%%%%%%%%%%%%%%%%%%%%%%%%%%%%%%%%%%%%%%%%%%%%%%%%%%%%%%%%%%%%%%%%%%%

\subsection{Standard}
\label{sec_finite_inters_1} Let $\cX_S = \cX_i$ for all $i \in
\Omega =[m]$, endowed with some probability $p_i \geq 0$. Since
$\cap_{i=1}^m \cX_i = \cX \subseteq \cX_S$,  then $\cX_S$ is a
stochastic approximation of $\cX$. Note that:
\[\cY=\left\{x\;:\; \sum_{i=1}^m p_i \Ib_{\cX_S}(x)=0\right\} = \bigcap_{i\;:\; p_i>0} \cX_i.\]
Hence,  a sufficient condition for exactness is to require  $p_i>0$
for all $i \in [m]$. Moreover, under this condition and
\eqref{linreg_finite} it follows that linear regularity
\eqref{linreg} holds with $\kappa =
\frac{\kappa_{\max}}{p_\text{min}}$, where $p_\text{min} = \min_{i
\in [m]} p_i$. Indeed, we can use the following inequality:
\[ p_\text{min} \max_{i \in [m]} \dist_{\cX_i}^2(x) \leq  \sum_{i=1}^m
p_\text{min} \dist_{\cX_i}^2(x) \leq  \sum_{i=1}^m p_i
\dist_{\cX_i}^2(x) = \Exp{\dist_{\cX_S}^2(x)}. \]

%%%%%%%%%%%%%%%%%%%%%%%%%%%%%%%%%%%%%%%%%%%%%%%%%%%%%%%%%%%%%%%

\subsection{Subsets}
\label{sec_finite_inters_2} With each nonempty subset $S \subseteq
[m] $ we associate a probability $p_S \geq 0$, such that $\sum_{S
\subseteq [m]} p_S = 1$. We then define  $\cX_S = \cap_{i\in S}
\cX_i$ with probability $p_S$. Since $\cX \subseteq \cX_S$, then
this is a stochastic approximation. Moreover,
\[\cY=\left\{x\;:\; \sum_{S} p_S \Ib_{\cX_S}(x) = 0 \right\} =\bigcap_{S\;:\;p_S>0} \cX_S.\]
A sufficient condition for the last set to be equal to $\cX$ (i.e.,
a sufficient condition for exactness) is \[[m]=\bigcup_{S\;:\;
p_S>0} S.\] In words, this condition requires us to assign positive
probabilities to some  collection of subsets covering $[m]$. If we
only assign positive probabilities to singletons,  we recover
example 1.  Moreover, under this condition and \eqref{linreg_finite}
it follows that linear regularity \eqref{linreg} holds with $\kappa
= \frac{\kappa_{\max}}{p_\text{min}}$, where $p_\text{min} =
\min_{S: p_S>0} p_S$. This is due to the fact that $\max_{i \in [m]}
\dist_{\cX_i}^2(x) \leq  \sum_{i=1}^m  \dist_{\cX_i}^2(x)$, that
$\dist_{\cX_i}^2(x) \leq \dist_{\cap_{j \in S} \cX_j}^2(x) =
\dist_{\cX_S}^2(x), \forall i \in S$ and that we assume  there is a
collection of subsets $S$ covering $[m]$.

%%%%%%%%%%%%%%%%%%%%%%%%%%%%%%%%%%%%%%%%%%%%%%%%%%%%%%%%%%%%%%%%%%%%

\subsection{Convex combination}
\label{sec_finite_inters_3} Fix $r \in [m]$, and  let us consider a
countable subset $\Omega_r$ defined as follows: \[  \Omega_r \subset
\left\{ S \in \R^{m}: \sum_{i=1}^m S_i =1, \; S \ge 0, \; \norm{S}_0
\le r \right\}. \] Let us consider a discrete probability
distribution $\Prob$ on $\Omega_r$.  We then choose $S \sim \Prob$
and define the stochastic approximation set as:
\[ \cX_S = \sum_{i=1}^m S_i\cX_i \eqdef \left\{ \sum_{i=1}^m S_i x_i \;:\; x_i \in \cX_i\right\}. \]
This is clearly a stochastic approximation, that is $\cX \subseteq
\cX_S$, since $\sum_{i=1}^m S_i =1$ and for any $x \in \cX$ it
follows that $x \in \cX_i$ for all $i \in [m]$ and thus $x = \sum_i
S_i x \in \cX_S$. For $r=1$ we recover the standard example from Section \ref{sec_finite_inters_1}. If additionally, we assume that  $\Omega_r$ contains
the basic vectors, i.e. $\{e_1, \cdots, e_m\} \subseteq \Omega_r$, and $\cX_S$ defined as
above, then exactness holds when $p_i = \Prob (S=e_i) > 0$ for all
$i \in [m]$. Indeed, if $x \in \cY$, then:
\[   0 = \Exp{ \Ib_{\cX_S}(x) } = \sum_{S \in \Omega}
p_S \Ib_{\cX_S}(x) \geq \sum_{S \in  \{e_1, \cdots, e_m\}} p_S
\Ib_{\cX_S}(x),    \] which implies $x \in \cX_i$, provided that
$p_i > 0$, for all $i \in [m]$. Moreover, under this condition and
\eqref{linreg_finite} it follows that linear regularity
\eqref{linreg} holds with $\kappa =
\frac{\kappa_{\max}}{p_\text{min}}$, where $p_\text{min} = \min_{i
\in [m]} p_i$. This is due to the fact that $\cX_{e_i} = \cX_i$ and
that:
\[ p_\text{min} \max_{i \in [m]} \dist_{\cX_i}^2(x) \leq
\sum_{i=1}^m p_i \dist_{\cX_i}^2(x) \leq \sum_{S \in \Omega} p_S
\dist_{\cX_S}^2(x) = \Exp{ \dist_{\cX_S}^2(x) }.
\]

%%%%%%%%%%%%%%%%%%%%%%%%%%%%%%%%%%%%%%%%%%%%%%%%%%%%%%%%%%%%%%%%%%%5

\subsection{Equality constraints}
\label{sec_finite_inters_4} Assume a linear representation for the
set $\cX$, that is  $\cX = \{x \in \rset^n: Ax = b  \}$, where the
matrix $A \in  \rset^{m \times n}$. In this case we have $\cX_i =
\{x \in \rset^n: A_i^T x - b_i = 0 \}$, where $A_i$ is the $i$th row
of matrix $A$. Let $q \le m$, $\Omega \subseteq   \rset^{m \times
q}$ and a probability distribution $\Prob$ on $\Omega$. Thus, we
define the stochastic approximation:
\[ \cX_S =\{x \in \rset^n: \;  S^T A x = S^T b \} \quad \forall S \in \Omega. \]
\textcolor{black}{We notice that  $\cap_{S \in \Omega} \cX_S = \{x:  S
A x = S b \; \forall S \in \Omega \}$. If we can find $m$ linearly
independent columns in the family of matrices $(S)_{S \in \Omega}$,
then $ \cX = \cap_{S \in \Omega} \cX_S $.} Next we derive sufficient
conditions for exactness, that is conditions that guarantee   $\cX=
\cY$, and we also provide an estimate for $\kappa$.

\begin{theorem}
\label{lemma_hyperplane} Let $\cX = \{x \in \rset^n: Ax=b\}$, with
$A \in \rset^{m \times n}$ and consider the stochastic approximation
$\cX_S = \{x \in \rset^n: S^T A x = S^T b \}$, where $S \in \rset^{m
\times q}$ is a random matrix in the probability space
$(\Omega,\Prob)$. Furthermore, assume that $S$ satisfies $\Exp{
S(S^TAA^TS)^{\dagger}S^T} \succ 0$. Then, we have exactness and the
linear regularity property \eqref{linreg} holds with constant:
\begin{equation}\label{kappa_linear}
\kappa = \frac{1}{\lambda_{\min}^{\text{nz}}(A^T\Exp{S(S^TA
A^TS)^{\dagger} S^T}A)} > 0.
\end{equation}
%where $\lambda_{\min}^{\text{nz}}(\cdot)$ denotes the smallest non-zero eigenvalue of a symmetric matrix.
\end{theorem}

%\begin{proof}
\noindent \textit{Proof}: Notice that the projection $\Pi_{\cX_S}(x)$ of $x$ onto $\cX_S$ can
be expressed as:
$$ \Pi_{\cX_S}(x) = x- A^T S(S^TAA^TS)^{\dagger} S^T(Ax-b),$$
thus the local distance $\dist_{\cX_S}(x)$ from $x$ to the set
$\cX_S$ is given by:
\begin{align}
\dist_{\cX_S}(x) = \norm{x - \Pi_{\cX_S}(x)}
&= \norm{A^TS(S^TAA^TS)^{\dagger} S^T(Ax-b)} \nonumber \\
&= \norm{A^TS(S^TAA^TS)^{\dagger} S^TA(x-\Pi_{\cX}(x))}.
\label{dist_xs}
\end{align}
Further, the matrix $P_S = A^TS(S^TA A^TS)^{\dagger} S^TA$ is
idempotent, that is $P_S^2= P_S$, which  implies that $\norm{P_S
z}^2 = z^T P_S z$ for any $z \in \rset^n$. By squaring and taking
expectation in both sides of \eqref{dist_xs} and also using the
previous property of $P_S$, we further obtain:
\begin{align}
\Exp{\dist_{\cX_S}^2(x)}
\overset{\eqref{dist_xs}}{=} & \Exp{\norm{P_S(x-\Pi_{\cX}(x))}^2  } \nonumber\\
= & \Exp{ (x-\Pi_{\cX}(x))^T P_S (x-\Pi_{\cX}(x)) } \nonumber\\
= & (x - \Pi_{\cX}(x))^T \Exp{ P_S }(x - \Pi_{\cX}(x)).
\label{interm_dist}
\end{align}
On the other hand, it is well known from the Courant-Fischer theorem
\cite{GolVan:96},  that for any $C  \in \rset^{m \times n}$ we have:
\begin{equation*}
\norm{C z} \ge \sigma^{\text{nz}}_{\min}(C)\norm{z}\qquad \forall z
\in \text{Im}(C^T),
\end{equation*}
where recall that  $\sigma^{\text{nz}}_{\min}$ denotes the smallest
nonzero singular value of a matrix. If we define the matrix $E =
\Exp{ S(S^TA A^TS)^{\dagger}S^T}$ and take $C = E^{1/2}A$, then the
above relation leads to:
\begin{equation}\label{bound_svalue}
\norm{E^{1/2} A z} \ge \sigma^{\text{nz}}_{\min}(E^{1/2}
A)\norm{z}\qquad \forall  z \in  \text{Im}(A^T E^{1/2}).
\end{equation}
Further, since we assume that $E= \Exp{S(S^TAA^TS)^{\dagger}S^T}
\succ 0$, then  $E^{1/2}  \succ 0$  and $\text{Im}(A^T) =
\text{Im}(A^T E^{1/2})$. Moreover,  we have the fact that $x -
\Pi_{\cX}(x) \in \text{Im}(A^T)$. Therefore, by applying the
relation \eqref{bound_svalue} for $z = x - \Pi_{\cX}(x)$, observing
that $\Exp{P_S} =A^T E A$, and by combining relations
\eqref{interm_dist} and \eqref{bound_svalue}, we have:
\begin{eqnarray*}
\Exp{\dist_{\cX_S}^2(x)}
&= & \norm{E^{1/2}A (x-\Pi_{\cX}(x))}^2  \\
&\overset{\eqref{bound_svalue}}{\ge}&
\left(\sigma_{\min}^{\text{nz}}(E^{1/2}A)\right)^2
\dist_{\cX}^2(x) \\
&=  & \lambda_{\min}^{\text{nz}}(A^T E A)\dist_{\cX}^2(x) \\
&=&  \lambda_{\min}^{\text{nz}} (\Exp{P_S})\dist_{\cX}^2(x) \\
&=  & \lambda_{\min}^{\text{nz}}(\Exp{A^TS(S^TA A^TS)^{\dagger}
S^TA})\dist_{\cX}^2(x)
\end{eqnarray*}
for all $ x \in \rset^n$. This final relation implies our statement. \qed
%\end{proof}

\noindent In \cite{GowRic:15} it has been proved that, when we
consider discrete samplings,  such as $S \in \Omega = \{e^1, \cdots,
e^m\}$, and full row rank matrices $A$ with no strictly zero rows,
the matrix $\Exp{S^T(SA A^TS^T)^{\dagger} S}$ is positive definite,
that is it satisfies our assumption considered in the previous
theorem. A simple consequence of previous theorem is the following:

\begin{corollary}
If we consider $\Omega = \{e_1, \cdots, e_m\}$, then for two choices
of the probability distribution $\Prob$ on $\Omega$  the linear
regularity constant takes the form:
\begin{align}
\label{const_lineq} \kappa  =
\begin{cases}
\frac{\norm{A}^2_F}{\lambda_{\min}^{\text{nz}}\left(A^TA\right)} &
\quad \text{if} \quad \Prob(S=e_i) = \frac{\norm{A_i}^2}{\norm{A}^2_F}\\
\frac{m}{\lambda_{\min}^{\text{nz}}\left(A^T D A\right)} & \quad
\text{if} \quad \Prob(S=e_i) = \frac{1}{m}
\end{cases}
\quad \geq 1,
\end{align}
where the diagonal matrix $D \eqdef \text{diag}(\|A_1\|^{-2},
\cdots, \|A_m\|^{-2})$.
\end{corollary}

%\begin{proof}
\noindent \textit{Proof}: If  $ \Omega = \{e_1, \cdots, e_m\}$ and the probability
 $\Prob(S=e_i) = \norm{A_i}^2/\norm{A}^2_F$, then the
stochastic approximation set $\cX_{e_i}$ is given by a linear
hyperplane, i.e. $\cX_{e_i} = \{x \in \rset^n: A_i^T x = b_i \}$,
and the expression in \eqref{kappa_linear} becomes:
\begin{align*}
& \lambda_{\min}^{\text{nz}}(A^T\Exp{S(S^TA A^TS)^{\dagger} S^T}A)
= \lambda_{\min}^{\text{nz}} \left( A^T\Exp{
\frac{e_ie_i^T}{\|A_i\|^2}}A \right) \\
& = \lambda_{\min}^{\text{nz}} \left( A^T \sum_{i=1}^m
\frac{\norm{A_i}^2}{\norm{A}^2_F} \frac{e_ie_i^T}{\|A_i\|^2}  A
\right) =  \lambda_{\min}^{\text{nz}} \left( A^T
\frac{I_m}{\|A\|_F^2} A \right) =
\frac{\lambda_{\min}^{\text{nz}}(A^TA)}{\|A\|^2_F}.
\end{align*}
Thus, in this case, the linear regularity constant is given by:
\begin{equation*}
\kappa \overset{\eqref{kappa_linear}}{=}
\frac{\norm{A}^2_F}{\lambda_{\min}^{\text{nz}}\left(A^TA\right)} =
\left(\frac{\norm{A}_F} {\sigma_{\min}^{\text{nz}}(A)}\right)^2 \geq
1.
\end{equation*}
For the uniform  probability  $\Prob(S=e_i) = 1/m$,  the expression
in \eqref{kappa_linear} becomes:
\begin{align*}
 \lambda_{\min}^{\text{nz}}(A^T\Exp{S(S^TA A^TS)^{\dagger} S^T}A)
& = \lambda_{\min}^{\text{nz}} \left( A^T\Exp{ \frac{e_ie_i^T}{\|A_i\|^2}}A \right) \\
& = \lambda_{\min}^{\text{nz}} \left( A^T \sum_{i=1}^m \frac{1}{m}
\frac{e_ie_i^T}{\|A_i\|^2}  A \right) =
\frac{\lambda_{\min}^{\text{nz}}(A^T D A)}{m},
\end{align*}
where the diagonal matrix $D = \text{diag}(\|A_1\|^{-2}, \cdots,
\|A_m\|^{-2})$. This proves our statement. \qed
% \end{proof}

%%%%%%%%%%%%%%%%%%%%%%%%%%%%%%%%%%%%%%%%%%%%%%%%%%%%%%%%%%%%%%%%%%%%%%%%%%%%%%%%%%%%

\subsection{Inequality constraints}
\label{sec_finite_inters_5} Let $q \le m$, $\Omega \subseteq
\rset^{m \times q}_+$ the set of matrices with nonnegative entries,
i.e., $\rset^{m \times q}_+ = \left\{ S \in \rset^{m \times q}:
S_{ij} \ge 0 \;\; \forall i \in [m], j \in [q] \right\}$, and a
probability distribution $\Prob$ on $\Omega$. Assume a functional
representation for the set $\cX$, that is  $\cX = \{x \in \rset^n:
\mathcal{F}(x) \le 0\}$, where $\mathcal{F}: \rset^n \to \rset^m$ is
a vector of convex closed functions, that is $\mathcal{F} =
(\mathcal{F}_1, \cdots, \mathcal{F}_m)$. In this case we have $\cX_i
= \{x \in \rset^n: \mathcal{F}_i(x) \le 0\}$. Thus, we define the
stochastic approximation:
\[ \cX_S =\{x \in \rset^n: S^T \mathcal{F}(x) \le 0\} \quad \forall S \in \Omega. \]
\textcolor{black}{ We notice that  $\cap_{S \in \Omega} \cX_S = \{x:
S^T \mathcal{F}(x) \le 0 \; \forall S \in \Omega \}$. If there exist
$m$ linearly independent columns in the family of matrices $(S)_{S
\in \Omega}$, then $ \cX = \cap_S \cX_S $. Moreover, if the
probability space is finite, then we also have exactness.} Next, we
provide estimates for the linear regularity constant $\kappa$ for
some particular sets. First, we consider finding a point in the
intersection of halfspaces, that is $\cX = \{x \in \rset^n: Ax \le b
\}$.

\begin{theorem}\label{lemma_halfspace_u}
Let $\cX = \{x \in \rset^n: Ax \le b \}$ and consider stochastic
approximation halfspaces $\cX_S =\{x \in \rset^n: S^TAx \le S^T b
\},$ where $S$ is a random vector from the finite probability space
$\Omega_r \subset \{S \in \rset^m: \; S \ge 0, \norm{S}_0 \le r \}$
for some given $r \in [m]$ endowed with a probability distribution
$\Prob =  (p_S)_{S \in \Omega_r}$. We further denote the Hoffman
constant for the polyhedral set $\cX$ with $\tilde{\kappa}$. Then,
under exactness the linear regularity property \eqref{linreg} holds
with constant:
\begin{equation}
\label{const_linineq_u} \kappa = \frac{\max\limits_{S \in \Omega_r}
\norm{A^T S}^2}{\min\limits_{S \in \Omega_r} \; p_S} \tilde{\kappa}.
\end{equation}
\end{theorem}

%\begin{proof}
\noindent \textit{Proof}: Notice that in this case we have an explicit projection onto $\cX_S$
given by $\Pi_{\cX_S}(x) = x -  \frac{\Pi_+(S^TAx -
S^Tb)}{\norm{A^TS}^2}A^TS$, which implies that:
\begin{equation}\label{bound_dist_lineq}
 \dist_{\cX_S}(x) = \frac{\Pi_+(S^TAx - S^Tb)}{\norm{A^TS}} \ge \frac{\Pi_+(S^TA x - S^Tb)}{
 \max\limits_{S \in \Omega_r} \norm{A^TS}}.
\end{equation}
From Markov inequality we have:
$$\frac{\Exp{\dist_{\cX_S}^2(x)}}{\max\limits_{S \in \Omega_r} \dist_{\cX_S}^2(x)}  \ge
 \Prob\left(\dist_{\cX_S}^2(x) \ge \max\limits_{S \in \Omega_r} \dist_{\cX_S}^2(x) \right).$$
Combining the previous inequality with \eqref{bound_dist_lineq}, we
obtain:
\begin{align}\label{bound_exp_lineq_u}
\Exp{\dist_{\cX_S}^2(x)} & \ge \Prob(\dist_{\cX_S}^2(x)
\ge\max\limits_{S \in \Omega_r} \dist_{\cX_S}^2(x))
\cdot \max\limits_{S \in \Omega_r} \dist_{\cX_S}^2(x) \nonumber\\
& = \Prob(\dist_{\cX_S}(x) \ge\max\limits_{S \in \Omega_r} \dist_{\cX_S}(x)) \cdot \max\limits_{S \in \Omega_r} \dist_{\cX_S}^2(x) \nonumber\\
& \ge \min\limits_{S \in \Omega_r} \; p_S \cdot \max\limits_{S \in \Omega_r} \dist_{\cX_S}^2(x) \nonumber  \\
& \overset{\eqref{bound_dist_lineq}}{\ge} \min\limits_{S \in
\Omega_r} \; p_S \cdot \frac{\max\limits_{S \in \Omega_r}
\Pi_+^2(S^TAx - S^Tb)}{\max\limits_{S \in \Omega_r} \norm{A^TS}^2}.
\end{align}

\noindent On the other hand it is well know that for a polyhedral
set the Hoffman inequality is valid, see \cite{BurFer:93}. Since we
assume exactness and that $\Omega_r$ has a finite number of
elements, then there exists some positive constant $\tilde{\kappa} >
0$ such that:
\[ \dist_{\cX}^2(x) \le \tilde{\kappa} \max\limits_{S \in \Omega_r} \Pi_+^2 (S^TAx - S^Tb)   \quad  \forall  x \in \rset^n.  \]
Using this Hoffman inequality  in \eqref{bound_exp_lineq_u} leads to
the relation:
\begin{equation*}
 \Exp{\dist_{\cX_S}^2(x)} \ge \frac{ \min\limits_{S \in \Omega_r} \; p_S }{\tilde{\kappa} \max\limits_{S \in \Omega_r}
  \norm{A^TS}^2}  \dist_{\cX}^2(x) \quad \forall x \in \rset^n,
\end{equation*}
which proves our statement. \qed
% \end{proof}

\noindent However, for a specific choice of the probability
distribution we can get better estimate for $\kappa$, as the next
corollary shows:
\begin{corollary}\label{lemma_halfspace}
Let $\cX = \{x \in \rset^n: Ax \le b \}$ and consider stochastic
approximation halfspaces $\cX_S =\{x \in \rset^n: S^TAx \le S^T b
\},$ where $S$ is a random vector from the finite probability space
$\Omega_r \subset \{S \in \rset^m: \; S \ge 0, \norm{S}_0 \le r \}$
for some given $r \in [m]$ endowed with the probability distribution
$\Prob =  (p_S)_{S \in \Omega_r}$ given by $p_S = \|A^TS\|^2/\sum_{S
\in \Omega_r} \|A^TS\|^2$. We further denote the Hoffman constant
for the polyhedral set $\cX$ with $\tilde{\kappa}$. Then, under
exactness the linear regularity property \eqref{linreg} holds with
constant:
\begin{equation}
\label{const_linineq} \kappa = \tilde{\kappa} \sum_{S \in \Omega_r}
\norm{A^T S}^2.
\end{equation}
\end{corollary}

%\begin{proof}
\noindent \textit{Proof}: Since $\Pi_{\cX_S}(x) = x -  \frac{\Pi_+(S^TAx -
S^Tb)}{\norm{A^TS}^2}A^TS$, then  we have:
\begin{equation*}
 \dist_{\cX_S}(x) = \frac{\Pi_+(S^TAx - S^Tb)}{\norm{A^TS}}.
\end{equation*}
Using the expressions for the distance and for the probability, we
further have:
\begin{align}\label{bound_exp_lineq}
\Exp{\dist_{\cX_S}^2(x)} & =  \sum_{S \in \Omega_r} p_S \dist_{\cX_S}^2(x)  \nonumber\\
& = \sum_{S \in \Omega_r} \frac{\|A^TS\|^2}{\sum_{S \in
\Omega_r}\|A^TS\|^2} \cdot \frac{\Pi_+^2(S^TAx -
S^Tb)}{\norm{A^TS}^2} \nonumber\\
&= \frac{1}{\sum_{S \in \Omega_r}\|A^TS\|^2} \sum_{S \in \Omega_r}
\Pi_+^2(S^TAx - S^Tb).
\end{align}

\noindent On the other hand, under exactness and  $\Omega_r$ has a
finite number of elements there exists some positive Hoffman
constant $\tilde{\kappa} > 0$ such that:
\[ \dist_{\cX}^2(x) \le \tilde{\kappa} \sum_{S \in \Omega_r} \Pi_+^2 (S^TAx - S^Tb)   \quad  \forall  x \in \rset^n.  \]
Using the Hoffman inequality  in \eqref{bound_exp_lineq} leads to
the relation:
\begin{equation*}
 \Exp{\dist_{\cX_S}^2(x)} \ge \frac{1}{\tilde{\kappa} \sum_{S \in \Omega_r}
  \norm{A^TS}^2}  \dist_{\cX}^2(x) \quad \forall x \in \rset^n,
\end{equation*}
which proves our statement. \qed
%\end{proof}

\noindent Second, following similar ideas as in \cite{Ned:10,GubPol:66}, we
consider the general case of a convex set $\cX$ with nonempty
interior,  that is, there exists a ball of radius $\delta > 0$ and
center  $\bar{x} \in \cX$ such that:
$$\{x \in \rset^n: \norm{\bar{x} - x} \le \delta \} \subseteq \cX.$$

\begin{theorem}\label{lemma_ball}
Let $\cX$ be a convex set with nonempty interior, that is there
exists $\delta > 0$ and $\bar{x} \in \cX$ such that $\{x \in
\rset^n: \norm{\bar{x} - x} \le \delta \} \subseteq \cX$.  Consider
any family of stochastic approximations $\cX_S$, where $S$ is a
random variable from the finite probability space $\Omega$ endowed
with a probability distribution $\Prob = (p_S)_{S \in \Omega}$.
Then, under exactness the linear regularity property \eqref{linreg}
holds over any bounded set $Q$ with constant:
\begin{equation}\label{const_ball}
\kappa = \frac{\max\limits_{x \in Q} \norm{x - \bar{x}}^2}{\delta^2
\min\limits_{S \in \Omega} \; p_S} \qquad \forall x \in Q.
\end{equation}
\end{theorem}

%\begin{proof}
\noindent \textit{Proof}: Let us define for some $\alpha > 0$ and $x \in \rset^n$  the vector:
$$y_{\alpha}(x) = \frac{\alpha}{\alpha + \delta} \bar{x} + \frac{\delta}{\alpha + \delta}  x.$$
Now we show that by choosing $\tilde{\alpha} = \max\limits_{S \in
\Omega} \dist_{\cX_S}(x)$,  then $y_{\tilde{\alpha}}(x) \in \cX$ for
all $x \in \rset^n$. Indeed, we first rewrite
$y_{\tilde{\alpha}}(x)$ as:
$$y_{\tilde{\alpha}}(x) = \frac{\tilde{\alpha}}{\tilde{\alpha} + \delta} z + \frac{\delta}{\tilde{\alpha}  + \delta} \Pi_{\cX_S}(x),$$
where $z = \bar{x} + \frac{\delta}{\tilde{\alpha}}(x -
\Pi_{\cX_S}(x))$. Notice that:
$$\norm{z - \bar{x}} = \frac{\delta}{\alpha}\norm{x - \Pi_{\cX_S}(x)} = \delta \frac{\dist_{\cX_S}(x)}{\max\limits_{S \in \Omega} \dist_{\cX_S}(x)} \le \delta.$$
Thus, we have $z \in \cX$, which implies $z \in \cX_S$ for all $S
\in \Omega$. Since $z \in \cX_S$, then further we conclude that also
$y_{\alpha}(x) \in \cX_S$ for all $\cX_S$, which finally confirms
that $y_{\alpha}(x) \in \cX$, due to exactness. By using this fact,
results that:
\begin{align}
\dist_{\cX}(x) &\le \norm{y_{\tilde{\alpha}}(x) - x}
 = \frac{\tilde{\alpha}}{\tilde{\alpha} + \delta} \norm{x - \bar{x}} \nonumber\\
& \le \frac{\tilde{\alpha}}{\delta} \norm{x - \bar{x}} =
\frac{\norm{x - \bar{x}}}{\delta} \max\limits_{S \in \Omega}
\dist_{\cX_S}(x). \label{max_errorbound}
\end{align}
From the Markov inequality we get the bound:
\begin{align}\label{markov}
\min\limits_{S \in \Omega} \; p_S \leq  \Prob(\dist_{\cX_S}^2(x) \ge
\max\limits_{S \in \Omega} \; \dist_{\cX_S}^2(x)) \le
\frac{\Exp{\dist^2_{\cX_S}(x)}}{\max\limits_{S \in \Omega}
\dist^2_{\cX_S}(x)}.
\end{align}
Using \eqref{max_errorbound} and \eqref{markov}, we obtain for any
$x \in Q$:
\begin{align*}
    \dist^2_{\cX}(x)
\overset{\eqref{max_errorbound} + \eqref{markov}} \le
\frac{\norm{x-\bar{x}}^2}{\delta^2 \min\limits_{S \in \Omega} \;
p_S}\Exp{\dist^2_{\cX_S}(x)} \le \frac{ \max\limits_{x \in Q}
\norm{x - \bar{x}}^2}{\delta^2 \min\limits_{S \in \Omega} \;
p_S}\Exp{\dist^2_{\cX_S}(x)},
\end{align*}
which confirms our result. \qed
%\end{proof}

%%%%%%%%%%%%%%%%%%%%%%%%%%%%%%%%%%%%%%%%%%%%%%%%%%%%%%%%%%%%%%%%%%%%%

%\subsection{Intersection of balls}
%\textcolor{black}{Assume $\cX_1,\dots,\cX_m$ are balls with varying
%centers and radii. Can we estimate $\kappa$.}

%%%%%%%%%%%%%%%%%%%%%%%%%%%%%%%%%%%%%%%%%%%%%%%%%%%%%%%%%%%%%%%%
%%%%%%%%%%%%%%%%%%%%%%%%%%%%%%%%%%%%%%%%%%%%%%%%%%%%%%%%%%%%%%%%

\section{Examples: infinite intersection}
Assume $\cX = \cap_{S\in
\Omega} \cX_S $, for some (possibly infinite) index set $\Omega$ and
sets $\cX_S\subseteq \R^n$. Many interesting applications can be modeled as the intersection of infinite (countable/uncountable) number of simple convex sets, see e.g. \cite{PatNec:17} for some control and machine learning applications. Let $\Prob$ be a probability measure on
$\Omega$. Then, if we choose $S\sim \Prob$,  $\cX_S$ is a stochastic
approximation of $\cX$. Note that
\[ \cY = \{x\;:\; \Prob(x\in \cX_S) = 1\}. \]
%For any $x\in \R^n$, let the event $\Omega(x) \eqdef \{S \in \Omega
%\;:\; x\in \cX_S\}$.

\subsection{Separation oracle}
Assume that we have access to a {\em
separation oracle} for $\cX$. That is, for each $S\in \R^n$, the
oracle either confirms that $S \in \cX$, or outputs a vector
$g=g(S)\in \R^n$ such that $\langle g, z- S \rangle \leq 0 $ for all
$z\in \cX$. If we let
\[\cX_S \eqdef \begin{cases} \R^n & \quad S\in \cX\\ \{x\;:\; \langle g, x- S \rangle \leq 0\} &
\quad  S\notin \cX, \end{cases}\] then clearly $\cX \subseteq \cX_S$
for all $S\in \R^n$. Given any distribution $\Prob$ over $\R^n$,
$\cX_S$ is a stochastic approximation of $\cX$.  \textcolor{black}{In
this case we can only  guarantee:
\[  \cX \subseteq \cap_{S \in \rset^n} \cX_S.  \]
}

%%%%%%%%%%%%%%%%%%%%%%%%%%%%%%%%%%%%%%%%%%%%%%%%%%%%%%%%%%%%%%

\subsection{Supporting halfspaces}
A particular case of the convex feasibility problem is the so-called
split feasibility problem \cite{ByrCen:01}: \[ \text{Find} \; x \in
\cX =\{x \in \rset^n: \; Ax \in \cZ \},
\]
i.e., $\cX$ is defined by imposing convex constraints defined by the
set $\cZ$ in the range of the matrix $A \in \rset^{m \times n}$.
Then, if we choose any $S \in \rset^n$ we can define a stochastic
approximation as the entire space or the following halfspace:
\[ \cX_S \eqdef
\begin{cases} \R^n & \quad S \in \cX\\ \{x:\;  c_S^T x \leq b_S\} & \quad  S \notin \cX, \end{cases}
\]
where $c_S \not =0$ and $b_S$ are defined as follows:
\[ c_S = A^T(AS - \Pi_\cZ(AS)) \; \text{and} \; b_S = \|AS\|^2 - (\Pi_\cZ(AS))^TAS - \|AS - \Pi_\cZ(AS)\|^2.
\]
Note that  the halfspace $\cX_S = \{x:\; c_S^T x \leq b_S\}$ can be
written equivalently as:
\[ \cX_S =  \{x:\; \langle AS - \Pi_\cZ(AS), Ax -
\Pi_\cZ(AS) \leq 0 \}.
\]
It is easy to check using the optimality conditions for the
projection onto $\cZ$ that for any $S \notin \cX$ the hyperplane
$c_S^T x = b_S$ separates $S$ from $\cX$, that is: \[ \cX \subseteq
\cX_S \quad \forall S \in \rset^n. \]  Therefore, given any
distribution $\Prob$ over $\R^n$, the halfspace $\cX_S$ is a
stochastic approximation of $\cX$. In fact, in this case we have:
\[  \cX = \cap_{S \in \rset^n} \cX_S.  \]
Indeed, it is straightforward that we have $\cX \subseteq \cap_{S
\in \rset^n} \cX_S$. For the other inclusion, let us take any $x \in
\cap_{S \in \rset^n} \cX_S$. Then, $x \in \cX_S$ for any fixed $S$.
Now, if we make the particular choice $S=x$, then $x \in \cX_x$,
that is it satisfies:
\[  \langle Ax - \Pi_\cZ(Ax), Ax -
\Pi_\cZ(Ax) \leq 0  \] which holds if and only if $Ax =
\Pi_\cZ(Ax)$, that is  $x \in \cX$. 

%\textcolor{black}{Can we
%generalize to sets $\cX = Y \cap Z$ and $Z$ is simple. }

%%%%%%%%%%%%%%%%%%%%%%%%%%%%%%%%%%%%%%%%%%%%%%%%%%%%%%%%%%%%%%%

\subsection{Normal cone} Let $\Omega\in \R^n$ be a closed convex set
and fix $\bar{x}\in \Omega$. Consider $\cX$ to be the normal cone of
the convex set $\Omega$ at some fixed point  $\bar{x} \in \Omega$:
\[\cX = \{x\;:\; (x-\bar{x})^T(S - \bar{x}) \leq 0 \; \text{for all} \; S\in \Omega\}
= \bigcap_{S\in \Omega} \cX_S,\] where $\cX_S \eqdef \{x\;:\;
(x-\bar{x})^T(S - \bar{x}) \leq 0\}$. If $\Prob$ is a probability
distribution over $\Omega$, and $S\sim \Prob$, then $\cX_S$ is a
stochastic approximation of $\cX$. Moreover, in this case we have
$\cX = \bigcap_{S\in \Omega} \cX_S$.

%\begin{remark}
%\noindent \textcolor{red}{We should add some comments here related
%to the fact that for future work we should investigate issues such
%as exactness and estimates for~$\kappa$. }
%\end{remark}

%%%%%%%%%%%%%%%%%%%%%%%%%%%%%%%%%%%%%%%%%%%%%%%%%%%%%%%%%%%%%
%%%%%%%%%%%%%%%%%%%%%%%%%%%%%%%%%%%%%%%%%%%%%%%%%%%%%%%%%%%

\section{Stochastic Projection Algorithm}
\label{sec_spa} In this section we propose the following parallel
stochastic projection method:

\vspace{10pt}

\begin{center}
\framebox{
\parbox{10.5 cm}{
\begin{center}
\textbf{ Algorithm SPA (general case) }
\end{center}
Choose $x^0 \in \R^n$, minibatch size $N \ge 1$, and positive
stepsizes $\{\alpha_k\}_{k \ge 0}$. For $k\geq 0$ repeat:
\begin{enumerate}
\item Draw $N$ independent samples, $S^k_1,S^k_2, \cdots, S^k_N \sim \Prob$
\item Compute $x^{k+1} = x^k - \alpha_k \left(x^k -  \frac{1}{N}\sum\limits_{i=1}^N
 \Pi_{\cX_{S^k_i}}(x^k)\right)$
\end{enumerate}
}}
\end{center}

\vspace{10pt}

This algorithm can be viewed as a random implementation of the
extrapolated method of parallel projections from \cite{Com:97},
which generates a sequence by  extrapolation of convex combinations
of  projections onto subfamilies of sets cyclically.

%%%%%%%%%%%%%%%%%%%%%%%%%%%%%%%%%%%%%%%%%%%%%%%%%%%%%%%%%%%%%%%%%%%%%%%

\subsection{Interpretation}
The minibatch algorithm \textbf{SPA} performs at each iteration $k$
a number of $N$ projections onto the simple sets $\cX_{S^k_1},
\cdots, \cX_{S^k_N}$ in parallel and then the new iterate is
computed taking a linear combination between the previous iterate
and the average of those projections. Such  minibatch strategy has
several interpretations. For example, when we consider the
stochastic smooth optimization problem \eqref{eq:reform_stoch_opt}:
\[ \min_{x \in \rset^n} F(x) = \Exp{F_S(x)}, \]
where $F_S(x) = 1/2 \| x - \Pi_{\cX_S}(x)\|^2$, usually  a Monte
Carlo simulation-based approach is used for solving it. It consists
in generating random samples of $S$ and the expected value function
$F$ is approximated by the corresponding sample average function.
That is, let $S_1, \cdots,S_N$ be independently and identically
distributed random sample of $N$ realizations of the random variable
$S$. Then, we consider the sample average function $\hat F_N = 1/N
\sum_{i=1}^N F_{S_i}$ and the associated problem:
\[ \min_{x \in \rset^n} \hat F_N(x). \]
Finally, this sample average optimization problem is solved. The
idea of using sample average approximations for solving stochastic
programs is a natural one and was used by various authors over the
years \cite{ShaDen:09}. However, the solution $\hat x_N^*$ of the
sample average optimization problem converges to the true solution
$x^*$ of the stochastic optimization problem only  for large enough
number of samples $N \to \infty$. On the other hand, in our
minibatch algorithm \textbf{SPA} the approach is different. First,
we fix the number of samples $N$. Then, at each iteration $k$ we
draw  $N$ independent samples $S^k_1,S^k_2, \cdots, S^k_N$ to also
form a sample average function $\hat F_N^k = 1/N \sum_{i=1}^N
F_{S_i^k}$. Finally, we do not solve the sample overage optimization
problem:
\[ \min_{x \in \rset^n} \hat F_N^k(x), \]
instead we only perform one gradient step for this problem with
stepsize $\alpha_k$
\[  x^{k+1} = x^k - \alpha_k \nabla \hat F_N^k(x^k),  \]
and then repeat the procedure. In this case we are not forced to
take $N$ large in order to obtain an approximative  solution of the
original problem. In fact, we can even consider $N=1$.

\noindent The minibatch algorithm \textbf{SPA} can be  also
interpreted in terms of the stochastic non-smooth optimization
problem \eqref{eq:reform_stoch_opt_indicator}:
\[ \min_{x \in \rset^n} f(x) = \Exp{f_S(x)}, \]
where $f_S(x) = \mathbb{I}_{\cX_{S}}(x)$. If we fix the number of
samples $N$, then at each iteration $k$ we draw  $N$ independent
samples $S^k_1,S^k_2, \cdots, S^k_N$ to form the same sample average
function:
\[ \hat F_N^k (x) = \frac{1}{N} \sum_{i=1}^N \left(\min_{z \in \rset^n} f_{S_i^k}(z) +
\frac{1}{2}\norm{x-z}^2 \right), \] and then consider solving the
sample overage optimization problem
\[ \min_{x \in \rset^n} \hat F_N^k(x), \]
which can be rewritten using the  notation $z = [z_1 \; \cdots \;
z_N]^T$ as follows:
\begin{equation*}
\min_{x \in \rset^n, z_i \in \rset^{n}} \; \hat F_N^k(x,z) \quad
\left( := \frac{1}{N} \sum\limits_{i=1}^N
\left[\mathbb{I}_{\cX_{S_i^k}}(z_i) +
\frac{1}{2}\norm{x-z_i}^2\right] \right).
\end{equation*}
However, we do not solve the previous average optimization problem
in the variables $(x,z)$, instead we only perform one  step of
Relaxed Block Alternating Minimization Method. That is,  given
$x^{k}$, we compute:
\begin{align*}
z^{k+1} = \arg\min\limits_{z \in \rset^{Nn}} & \hat  F_N^k
(x^{k},z), \quad \tilde{x}^{k+1} = \arg\min\limits_{x \in \rset^{n}} \hat F_N^k(x,z^{k+1}) \\
x^{k+1} &= (1 - \alpha_k) x^k + \alpha_k \tilde{x}^{k+1},
\end{align*}
and repeat the whole procedure. Again, this strategy allows us to
work also with $N$ small, including  $N=1$.

%%%%%%%%%%%%%%%%%%%%%%%%%%%%%%%%%%%%%%%%%%%%%%%%%%%%%%%%%%%%%%%%%%%%%%%%

\subsection{Convergence analysis}
Our convergence analysis is based on two important properties of the
family of convex sets $(\cX_S)_{S \in \Omega}$. For simplicity, we
recall them once more here. First, there exists $\Uconst \leq 1$
satisfying the inequality \eqref{eq:gamma0}, i.e.:
\begin{equation}
\label{eq:gamma} \left\|\Exp{x-\Pi_{\cX_{S}}(x)}\right\|^2 \leq
\Uconst \cdot \Exp{ \left\| x - \Pi_{\cX_{S}}(x) \right\|^2 } \qquad
\forall x\in \R^n.
\end{equation}
However, for specific sets and distributions $\Prob$, we  proved in
Section \ref{sec_F} that   $\Uconst$ can be much smaller than $1$.
Second, there exists $\kappa \leq \infty$ such that the family of
convex sets $(\cX_S)_{S \in \Omega}$ satisfies the linear regularity
property~\eqref{linreg}, i.e.:
\begin{equation}
\label{linreg0} \dist_{\cX}^2(x) \le \kappa \; \Exp{
\dist_{\cX_S}^2(x)}  \qquad \forall x\in \R^n.
\end{equation}
However, we have proved in Section \ref{sec_finite_inters} that for
specific sets and distributions $\Prob$,  the constant  $\kappa$ can
be finite, that is   $\kappa < \infty$. Based on the properties
\eqref{eq:gamma} and \eqref{linreg0} the smooth objective function
$F$ of the stochastic optimization problem
\eqref{eq:reform_stoch_opt}  satisfies Theorem \ref{th_shapeF}, in
particular we have:
\begin{align}
\label{qg_Lip3} \frac{1}{2 \kappa} \| x - \Pi_{\cX}(x) \|^2 \leq
F(x) -F^* \leq \frac{\Uconst}{2} \| x - \Pi_{\cX}(x) \|^2 \qquad
\forall x \in \rset^n.
\end{align}
There is an interesting interpretation of  inequality
\eqref{qg_Lip3}, that is  the objective function $F$ is strongly
convex with constant $\frac{1}{\kappa}$ and has Lipschitz continuous
gradient with constant $\Uconst \leq 1$ when restricted along any
segment $[x, \Pi_{\cX}(x)]$. Thus,  $\kappa \Uconst$ represents  the
condition number of the convex feasibility problem
\eqref{convexfeas_original}. Using the inequalities
\eqref{eq:gamma}-\eqref{qg_Lip3} we can prove not only asymptotic
convergence of the sequence $\{ x^k \}_{k \ge 0}$ generated by
algorithm \textbf{SPA}, but also rates of convergence. We start with
a basic result from probability theory, see e.g. \cite{Ned:11}:

\begin{lemma}
[Supermartingale Convergence Lemma]  Let $v^k$ and $u^k$ be
sequences of nonnegative random variables such that:
\[ \Exp{v^{k+1} | F_k} \leq v^k - u^k  \quad \text{a.s.} \quad \forall k \geq 0,  \]
where $F_k$ denotes the collection $\{v^0, \cdots, v^k, u^0, \cdots,
u^k\}$. Then, we have $v^k$ convergent to a random variable $v$ a.s.
and $\sum_{k=0}^\infty u^k < \infty$ a.s.
\end{lemma}

\noindent Then, we obtain the following asymptotic convergence
result:

\begin{theorem}
\label{thm:convergence} Assume that the set $\cX$ is nonempty and
define $\Uconst_N \eqdef \frac{1}{N} +
\left(1-\frac{1}{N}\right)\Uconst \leq 1$. Let $\{ x^k \}_{k \ge 0}$
be generated by algorithm \textbf{SPA} with stepsizes $0 < \alpha_k
< \frac{2}{\Uconst_N}$. Then, we have the following average
decrease:
\begin{align}
\label{central_ineq} \Exp{\|x^{k+1}-x^*\|^2 \;|\; x^k} \leq
\|x^k-x^*\|^2 - 2 (2\alpha_k - \alpha_k^2 \Uconst_N) F(x^k)
\end{align}
for all $k\geq 0$ and $x^* \in \cX$.  Moreover, the fastest decrease
is given by the constant stepsize $\alpha_k=1/\Uconst_N$. If
additionally, exactness holds and the stepsize satisfies $\delta
\leq \alpha_k \leq \frac{2}{\Uconst_N} - \delta$ for some $0< \delta
\leq \frac{1}{\Uconst_N}$, then the sequence $x^k$ converges almost
sure to a random point in the set $\cX$ and $\lim\limits_{k \to
\infty} F(x^k) =0$ almost sure.
\end{theorem}

%\begin{proof}
\noindent \textit{Proof}: For simplicity, we shall write $\Pi^k_i=\Pi_{\cX_{S^k_i}}(x^k)$. Let
$x^*$ be any element of $\cX$. Then, we have the following:
\begin{align}
\label{eq:decrease} \|x^{k+1} \!- x^*\|^2 \! & \!= \left\|x^k - x^*
- \alpha_k \left( x^k - \frac{1}{N} \sum_{i=1}^N \Pi^k_i
\right)\right\|^2 \nonumber \\
& =\left\|x^k - x^* - \alpha_k \frac{1}{N} \sum_{i=1}^N (x^k - \Pi^k_i )\right\|^2 \nonumber\\
&= \! \|x^k -x^*\|^2 -\! \frac{2\alpha_k}{N} \! \sum_{i=1}^N \!
\left\langle x^k \!- x^* , x^k \!- \Pi^k_i \right
\rangle \!+\! \frac{\alpha_k^2}{N^2} \! \left\|\sum_{i=1}^N  (x^k \!- \Pi^k_i) \right\|^2 \nonumber\\
&\leq  \|x^k -x^*\|^2 - \frac{2\alpha_k}{N}\sum_{i=1}^N  \left\|x^k
- \Pi^k_i \right\|^2 +
\frac{\alpha_k^2}{N^2} \left\|\sum_{i=1}^N  (x^k - \Pi^k_i) \right\|^2 \nonumber\\
&= \|x^k -x^*\|^2 - \frac{2\alpha_k}{N}\sum_{i=1}^N  \left\|x^k -
\Pi^k_i \right\|^2  \\
& \qquad \qquad +   \frac{\alpha_k^2}{N^2} \left(\sum_{i=1}^N \|x^k
- \Pi^k_i\|^2 + \sum_{i\neq j} \langle x^k - \Pi^k_i, x^k - \Pi^k_j
\rangle \right), \nonumber
\end{align}
where the inequality follows from the bound: \[ \langle x^k - x^*,
x^k - \Pi_i^k  \rangle = \langle x^k - \Pi_i^k, x^k - \Pi_i^k
\rangle + \langle \Pi_i^k - x^*, x^k - \Pi_i^k   \rangle \geq \| x^k
- \Pi_i^k \|^2, \] since $ \langle \Pi_i^k - x^*, x^k - \Pi_i^k
\rangle \geq 0$ for all $x^* \in \cX \subseteq \cX_{S_i^k}$. Taking
expectations conditioned on $x^k$ and using the definition of $F$:
\begin{align*}
F(x^k) & = \frac{1}{2} \Exp{\|x^k- \Pi^k_i\|^2\;|\; x^k} =
\frac{1}{2} \Exp{\|x^k - \Pi_{\cX_{S^k_i}}(x^k)\|^2\;|\; x^k} \\
& = \frac{1}{2} \Exp{\|x^k- \Pi_{\cX_S}(x^k)\|^2\;|\; x^k},
\end{align*}
and invoking conditional independence of $\Pi^k_i$ and $\Pi^k_j$ for
$i \neq j$ (inherited from independence of $S^k_i$ and $S^k_j$), we
obtain:
\begin{align}
&\Exp{ \|x^{k+1} - x^*\|^2  \;|\; x^k}
 \leq \|x^k-x^*\|^2 - 4 \alpha_k F(x^k) \nonumber \\
& \qquad \qquad + \frac{\alpha_k^2}{N^2} \left(2 N F(x^k) + \sum_{i\neq j}
\left\langle \Exp{x^k - \Pi^k_i\;|\; x^k}, \Exp{x^k - \Pi^k_j\;|\; x^k} \right\rangle \right) \nonumber \\
&= \|x^k \!- x^*\|^2 \!- 4 \alpha_k F(x^k) \!+ \frac{2\alpha_k^2}{N}F(x^k) \!+\frac{\alpha_k^2 (N^2\!-\!N)}{N^2} \left\| \Exp{x^k \!-\! \Pi_{\cX_S}(x^k)\;| x^k} \right\|^2
\label{descent_prelim}\\
&\overset{\eqref{eq:gamma}}{\leq} \|x^k \!-x^*\|^2 \!- 4 \alpha_k
F(x^k) \!+ \frac{2 \alpha_k^2}{N}F(x^k) \!+ \frac{\alpha_k^2 (N\!-\!1)}{N} \Uconst \Exp{\| x^k \!-\! \Pi_{\cX_S}(x^k) \|^2 \;| x^k} \nonumber \\
&=  \|x^k-x^*\|^2 - 4 \alpha_k F(x^k)
+ \frac{2\alpha_k^2}{N}F(x^k) + \frac{2\alpha_k^2 (N-1)}{N} \Uconst F(x^k) \nonumber \\
&= \|x^k-x^*\|^2 - 2 (2\alpha_k - \alpha_k^2 \Uconst_N) F(x^k). \nonumber
\end{align}
Thus, we have obtained  for all $k \geq 0$ and  $x^* \in \cX$:
\begin{align*}
\Exp{\|x^{k+1}-x^*\|^2 \;|\; x^k} \leq \|x^k-x^*\|^2 - 2 (2\alpha_k
- \alpha_k^2 \Uconst_N) F(x^k).
\end{align*}
Clearly, the fastest  decrease is obtained by maximizing $2\alpha_k
- \alpha_k^2 \Uconst_N$ in $\alpha_k$, that is the maximum  is
obtained for constant stepsize $\alpha_k=1/\Uconst_N$. Further, for the
stepsizes  satisfying $\delta \leq \alpha_k \leq \frac{2}{\Uconst_N} - \delta$  we have
$2  \alpha_k - \alpha^2_k \Uconst_N \geq  \delta^2 \Uconst_N >  0$.
Then, from  Supermartingale Convergence Lemma we have that $\|x^k-x^*\|^2$
converges a.s. for every $x^* \in \cX$ and thus the sequence $x^k$
is bounded a.s. This implies that $x^k$ has a limit point $\tilde
x^*$.  Since we also have $\sum_{k=0}^\infty F(x^k) < \infty$ a.s.,
it follows that $F(x^k) \to 0$ a.s. Therefore, for any accumulation
point $\tilde x^*$ of $x^k$ we have $F(\tilde x^*) =0$ a.s. (by
continuity of $F$). This leads to $\tilde x^* \in \cY$ a.s.  When
exactness holds (i.e. $\cX=\cY$), it follows that at least a
subsequence of $x^k$ converges almost surely to a random point
$\tilde x^*$ from the set $\cX$. \qed
%\end{proof}

\noindent The previous theorem clearly shows that in order to have
decrease in average distances (see \eqref{central_ineq})
the stepsize $\alpha_k$ has to satisfy:
\begin{align}
\label{eq:stepsize} 0 < \alpha_k  < \frac{2}{\Uconst_N} \qquad
\forall k \geq 0.
\end{align}
This shows that we can use large stepsizes $\alpha_k$. Thus, we
prove theoretically, what is known in numerical applications for a
long time, namely that this overrelaxation $\alpha_k  \approx
\frac{2}{\Uconst_N} > 1$ accelerates significantly in practice the
convergence of projection methods as compared to its basic
counterpart $\alpha_k=1$, see \cite{CenChe:12,Com:97}. For several
important sets we can  estimate the ``Lipschitz'' constant
$\Uconst$ and consequently $\Uconst_N$, see Section \ref{sec_F}. For
other sets however, it is difficult to compute $\Uconst$. In this
case we propose an adaptive estimation of  $\Uconst$ at each
iteration $k$ as follows:
\[  \Uconst^k =  \frac{\|\Exp{x^k -\Pi_{\cX_{S}}(x^k)}\|^2}{
\Exp{\| x^k - \Pi_{\cX_{S}}(x^k)\|^2}}. \] This choice has the
following interpretation. From Theorem \ref{th_shapeF} we have that
$F$ has Lipschitz continuous gradient with constant $\Uconst$ on any segment
$[x^k, \ \Pi_{\cX}(x^k)]$: \[ F(x^k) \overset{\eqref{qg_Lip2}}{\geq} F^* + \langle \nabla F(\Pi_{\cX}(x)), x^k - \Pi_{\cX}(x)
\rangle + \frac{1}{2 \Uconst} \| \nabla F(x^k) - \nabla F(\Pi_{\cX}(x))\|^2, \]
which, using $F^* = F(\Pi_{\cX}(x^k)) = 0$ and $\nabla
F(\Pi_{\cX}(x^k)) = 0$, is equivalent to:
\[ \Uconst \geq \frac{ 1/2 \| \nabla F(x^k)\|^2 }{F(x^k)} = \Uconst^k.   \]
Using arguments of Theorem \ref{thm:convergence}, it is straightforward to obtain the following descent.

\begin{corollary}
\label{cor:convergence} Assume that the set $\cX$ is nonempty and
define $\Uconst_N^k \eqdef \frac{1}{N} +
\left(1-\frac{1}{N}\right)\Uconst^k \leq 1$. Let $\{ x^k \}_{k \ge 0}$
be generated by algorithm \textbf{SPA} with stepsizes $0 < \alpha_k
< \frac{2}{\Uconst_N^k}$. Then, we have the following average
decrease:
\begin{align}
\label{central_ineq} \Exp{\|x^{k+1}-x^*\|^2 \;|\; x^k} \leq
\|x^k-x^*\|^2 - 2 (2\alpha_k - \alpha_k^2 \Uconst_N^k) F(x^k)
\end{align}
for all $k\geq 0$ and $x^* \in \cX$.
\end{corollary}

%\begin{proof}
\noindent \textit{Proof}: From the relation \eqref{descent_prelim} we have:
\begin{align*}
\Exp{ \|x^{k+1} - x^*\|^2  \;|\; x^k}
& \leq  \|x^k \!- x^*\|^2 \!- 4 \alpha_k F(x^k) \!+ \frac{2\alpha_k^2}{N}F(x^k) \!+\frac{\alpha_k^2 (N^2\!-\!N)}{N^2} \left\| \Exp{x^k \!-\! \Pi_{\cX_S}(x^k)\;| x^k} \right\|^2\\
& = \|x^k \!-x^*\|^2 \!- 4 \alpha_k
F(x^k) \!+ \frac{2 \alpha_k^2}{N}F(x^k) \!+ \frac{\alpha_k^2 (N\!-\!1)}{N} \gamma^k \Exp{\| x^k \!-\! \Pi_{\cX_S}(x^k) \|^2 \;| x^k} \nonumber \\
&=  \|x^k-x^*\|^2 - 4 \alpha_k F(x^k) + \frac{2\alpha_k^2}{N}F(x^k) + \frac{2\alpha_k^2 (N-1)}{N} \gamma^k F(x^k) \nonumber \\
&= \|x^k-x^*\|^2 - 2 (2\alpha_k - \alpha_k^2 \gamma_N^k) F(x^k), \nonumber
\end{align*}
which confirms the results.
\qed

% Thus, we have obtained  for all $k \geq 0$ and  $x^* \in \cX$:
% \begin{align*}
% \Exp{\|x^{k+1}-x^*\|^2 \;|\; x^k} \leq \|x^k-x^*\|^2 - 2 (2\alpha_k
% - \alpha_k^2 \Uconst_N^k) F(x^k).
% \end{align*}
% Clearly, the fastest  decrease is obtained by maximizing $2\alpha_k
% - \alpha_k^2 \Uconst_N^k$ in $\alpha_k$, that is the maximum  is
% obtained for constant stepsize $\alpha_k=1/\Uconst_N^k$. Further, for the
% stepsizes  satisfying $\delta \leq \alpha_k \leq \frac{2}{\Uconst_N^k} - \delta$  we have
% $2  \alpha_k - 	\alpha^2_k \Uconst_N \geq  \delta^2 \Uconst_N^k >  0$.\qed
%\end{proof}

When even the expectation is difficult to compute for finding
$\Uconst^k$, then, inspired by \cite{Com:97}, we propose to use the following approximation
for the previous ratio:
\[  \Uconst^k = \frac{\|  \sum_{i=1}^N w_i^k (x^k -\Pi_{\cX_{S_i^k}}(x^k)) \|^2}{
\sum_{i=1}^N w_i^k  \| x^k -  \Pi_{\cX_{S_i^k}}(x^k) \|^2},
\]
where the weights $w_i^k$ satisfy $ \sum_{i=1}^N w_i^k=1$ and $w_i^k
>0$. For these situations we %, we define $\Uconst^k_N= \frac{1}{N} + (1 - \frac{1}{N}) \Uconst^k$ and
take the stepsize:
\[  \alpha_k = \frac{\alpha}{\Uconst^k_N}, \quad \text{with} \; \alpha \in (0, \ 2).   \]
The effectiveness of this choice for the stepsize has been shown in
many practical applications, see e.g. \cite{CenChe:12,Com:97}.  Next
theorem provides rates of convergence for the sequence $x^k$
generated by \textbf{SPA}:

\begin{theorem} \label{thm:general}
Assume that the set $\cX$  is nonempty and
define $\Uconst_N = \frac{1}{N} +\left(1-\frac{1}{N}\right)\Uconst$.
Let $\{ x^k \}_{k \ge 0}$ be generated by algorithm \textbf{SPA}
with stepsizes satisfying $\delta \leq \alpha_k \leq
\frac{2}{\Uconst_N} - \delta$ for some $0< \delta \leq
\frac{1}{\Uconst_N}$. Then:
\begin{enumerate}
\item[$(i)$] For the average point $\hat{x}^k = \frac{1}{\Sigma_k}\sum\limits_{i=0}^{k-1} \alpha_i
x^i$, where $\Sigma_k = \sum\limits_{i=0}^{k-1} \alpha_i$, we have
the following sublinear convergence rate:
\begin{equation*}
\Exp{F(\hat{x}^k)} - F^* =  \frac{1}{2} \Exp{
\dist^2_{\cX_{S}}(\hat{x}^k)} \leq \frac{ \dist^2_{\cX}(x^{0})}{2
\delta \Uconst_N  \Sigma_k}.
\end{equation*}
Moreover,  the average sequence $\hat{x}^k$ converges almost surely
to a random point in the set~$\cX$, provided that exactness holds.

\item[$(ii)$] If additionally  the linear regularity property \eqref{linreg0}
holds, then we have the following linear convergence  rate for the
last iterate $x^k$:
\begin{equation*}
\Exp{\dist^2_{\cX}(x^{k+1})} \leq \left(1- \frac{\delta^2
\Uconst_N}{ \kappa } \right) \Exp{\dist^2_{\cX}(x^{k})},
\end{equation*}
or in terms of function values:
\begin{equation*}
\Exp{F(x^{k})} - F^* \leq    \left(1- \frac{\delta^2
\Uconst_N}{\kappa } \right)^k \frac{\Uconst \dist^2_{\cX}(x^0)}{2}.
\end{equation*}
\end{enumerate}
\end{theorem}

%\begin{proof}
\noindent \textit{Proof}: By taking expectation w.r.t. the entire history on both sides in \eqref{central_ineq} we get the following decrease in the distance
to a point  $x^* \in \cX$:
\begin{align*}
\Exp{ \|x^{k+1}-x^*\|^2 } \leq \Exp{ \|x^k-x^*\|^2 } - 2(2\alpha_k -
\alpha_k^2 \Uconst_N) \Exp{ F(x^k)}.
\end{align*}
Further,  denoting $r_k \eqdef \Exp{\norm{x^k - x^*}^2}$ and
noticing the lower bound $2 - \alpha_k \Uconst_N \geq \delta
\Uconst_N$ for any stepsize satisfying $\delta \leq \alpha_k \leq
\frac{2}{\Uconst_N} - \delta$ for some $0< \delta \leq
\frac{1}{\Uconst_N}$, we have:
\begin{equation*}
2 \delta \Uconst_N \alpha_k \Exp{F(x^k)} \leq 2 \alpha_k (2 -
\alpha_k \Uconst_N) \Exp{F(x^k)} \le r_k - r_{k+1}.
\end{equation*}
If we add the entire history from $i=0$ to $i= k-1$, we obtain:
\begin{equation*}
2 \delta \Uconst_N  \Exp{\sum\limits_{i=0}^{k-1} \alpha_i F(x^i)}
\!=\! \sum\limits_{i=0}^{k-1} 2 \delta \Uconst_N \alpha_i
\Exp{F(x^i)} \!\le r_0 - r_{k} \!\le r_0 =\! \norm{x^0 -x^*}^2
\end{equation*}
for all $x^* \in \cX$. If we choose $x^* = \Pi_{\cX}(x^0)$ and use
the  convexity of function $F$, then we finally get:
\begin{equation*}
2 \Sigma_k \delta \Uconst_N  \Exp{
F\left(\frac{1}{\Sigma_k}\sum\limits_{i=0}^{k-1} \alpha_i
x^i\right)} \le 2 \delta \Uconst_N \Exp{ \sum\limits_{i=0}^{k-1}
\alpha_i F(x^i) } \le \dist_{\cX}^2(x^0).
\end{equation*}
This relation and $F^*=0$ imply  immediately the first part of our
result. Moreover,  by Theorem \ref{thm:convergence}, $x^k$ converges
almost surely to a random point in the set $\cX$. Therefore, the
average sequence  $\hat{x}^k = \frac{1}{k}\sum\limits_{i=0}^{k-1}
x^i$ also converges almost surely to the same random point in the
set $\cX$.

\noindent (ii) In order to prove linear convergence under linear
regularity property \eqref{linreg0} we use again inequality
\eqref{central_ineq} and $F^*=0$:
\begin{eqnarray*}
\Exp{\|x^{k+1}-x^*\|^2 \;|\; x^k}
 \leq & \|x^k-x^*\|^2 - 2 (2\alpha_k - \alpha_k^2 \Uconst_N) \left( F(x^k) - F^* \right) \\
 \overset{\eqref{qg_Lip3}}{\leq} & \hspace{-1.5cm}\|x^k-x^*\|^2 -  \frac{2\alpha_k - \alpha_k^2 \Uconst_N}{\kappa} \dist^2_{\cX}(x^k).
\end{eqnarray*}
Taking expectations w.r.t. the entire history, we obtain:
\begin{equation}\label{eq:9ys09u09}
\Exp{\|x^{k+1}-x^*\|^2} \leq  \Exp{\|x^k-x^*\|^2} - \frac{2\alpha_k
- \alpha_k^2 \Uconst_N}{\kappa} \Exp{\dist^2_{\cX}(x^k)}.
\end{equation}
Choosing $x^* = \Pi_{\cX}(x^k)$, and using the inequality
$\dist^2_{\cX}(x^{k+1}) = \|x^{k+1}-\Pi_{\cX}(x^{k+1})\|^2 \leq
\|x^{k+1}-x^*\|^2$ together with \eqref{eq:9ys09u09},  we finally
get:
\[\Exp{\dist^2_{\cX}(x^{k+1})} \leq \left(1- \frac{2\alpha_k - \alpha_k^2 \Uconst_N}{\kappa}
\right) \Exp{\dist^2_{\cX}(x^{k})}.\] Since for our choice of the
stepsize  $\delta \leq \alpha_k \leq \frac{2}{\Uconst_N} - \delta$
for some $0< \delta \leq \frac{1}{\Uconst_N}$, we have $2\alpha_k -
\alpha_k^2 \Uconst_N \geq \delta^2 \Uconst_N$, then  the previous
relation  implies immediately:
\[\Exp{\dist^2_{\cX}(x^{k+1})} \leq \left(1- \frac{\delta^2 \Uconst_N}{\kappa}
\right) \Exp{\dist^2_{\cX}(x^{k})}.\] which proves the second
statement of the theorem. Finally, combining the convergence rate in
distances with the right hand side inequality in \eqref{qg_Lip3} we
get the convergence in expectation of value function.  \qed
%\end{proof}

\noindent An immediate consequence of Theorem \ref{thm:general} is
the following corollary:
\begin{corollary}
\label{thm:general_opt} Assume that the set $\cX$  is nonempty and
$\Uconst_N = \frac{1}{N} + \left(1-\frac{1}{N}\right)\Uconst$. Let
$\{ x^k \}_{k \ge 0}$ be generated by algorithm \textbf{SPA} with
the optimal constant stepsize $\alpha_k=1/\Uconst_N$. Then:
\begin{enumerate}
\item[$(i)$]  For the average point
$\hat{x}^k = \frac{1}{k}\sum\limits_{i=0}^{k-1} x^i$ we have the
following sublinear convergence rate:
\begin{equation*}
\Exp{F(\hat{x}^k)} - F^* =  \frac{1}{2} \Exp{
\dist^2_{\cX_{S}}(\hat{x}^k)} \leq \frac{\Uconst_N \cdot
\dist^2_{\cX}(x^{0})}{2 k}.
\end{equation*}

\item[$(ii)$] If additionally  the linear regularity property \eqref{linreg0}
holds, then  we have the following linear convergence  rate for the
last iterate $x^k$:
\begin{equation}
\label{SPA_linconv} \Exp{\dist^2_{\cX}(x^{k+1})} \leq \left(1-
\frac{1}{\Uconst_N \cdot \kappa } \right)
\Exp{\dist^2_{\cX}(x^{k})},
\end{equation}
or in terms of function values:
\begin{equation*}
\Exp{F(x^{k})} - F^* \leq    \left(1- \frac{1}{\Uconst_N \cdot
\kappa } \right)^k \frac{\Uconst \dist^2_{\cX}(x^0)}{2}.
\end{equation*}
\end{enumerate}
\end{corollary}

%\begin{proof}
\noindent \textit{Proof}: By taking expectation w.r.t. the entire history on both sides in
\eqref{central_ineq} we get the following decrease in the distance
to a point  $x^* \in \cX$:
\begin{align*}
%\label{central_ineq}
\Exp{ \|x^{k+1}-x^*\|^2 } \leq \Exp{ \|x^k-x^*\|^2 } - 2(2\alpha_k -
\alpha_k^2 \Uconst_N) \Exp{ F(x^k)}.
\end{align*}
Further,  denoting $r_k \eqdef \Exp{\norm{x^k - x^*}^2}$, we have:
\begin{equation*}
2 (2\alpha_k - \alpha_k^2 \Uconst_N) \Exp{F(x^k)} \le r_k - r_{k+1}.
\end{equation*}
The fastest decrease is obtained maximizing $2\alpha_k - \alpha_k^2
\Uconst_N$ in $\alpha_k$, which leads to the optimal stepsize
$\alpha_k = 1/\Uconst_N$. The rest of the proof follows exactly the
same steps as in the proof of Theorem \ref{thm:general}, observing
that choosing $\delta = 1/\Uconst_N$ we get $\alpha_k =
1/\Uconst_N$. \qed

\noindent From Theorem \ref{thm:general} and Corollary
\ref{thm:general_opt} it follows that the convergence rates of
algorithm \textbf{SPA} depend explicitly on the minibatch sample
size  $N$ via the term $\Uconst_N$. Moreover, we notice that the
scheme \textbf{SPA} is very  general and  we can  recover multiple
existing projection algorithms from the literature. Further, we
analyze some particular algorithms  resulted from \textbf{SPA} and
derive  their convergence~rates.

%%%%%%%%%%%%%%%%%%%%%%%%%%%%%%%%%%%%%%%%%%%%%%%%%%%%%%%%%%%%%%%%%%

\subsection{Average Projection algorithm: $N = m/\infty$}
\noindent  As $N \to  \infty$, we have
$\Uconst_N \to \Uconst$, and the linear rate in Corollary
\ref{thm:general_opt}  converges to $1-1/(\kappa \cdot \Uconst)$.
This is also confirmed by the convergence rate given in Theorem
\ref{th:conv_GA} below. More precisely, when $N \to \infty$ the
algorithm \textbf{SPA} becomes the deterministic gradient method for
solving the smooth convex problem \eqref{eq:reform_stoch_opt},
which we call average projection algorithm:

\vspace{10pt}

\begin{center}
\framebox{
\parbox{12 cm}{
\begin{center}
\textbf{ Algorithm AvP  }
\end{center}
Choose $x^0 \in \R^n$ and positive stepsizes $\{\alpha_k\}_{k \ge 0}$. For
$k\geq 0$ repeat:
\begin{enumerate}
\item Compute $x^{k+1} = x^k - \alpha_k \nabla F(x^k) \quad \left(\;\; \eqdef x^k - \alpha_k
\left( x^k - \Exp{\Pi_{\cX_S}(x^k)} \right)\right)$
\end{enumerate}
}}
\end{center}

\vspace{10pt}

\noindent Under linear regularity condition \eqref{linreg0} the
sequence $\{x^k\}_{k \ge 0}$ generated by algorithm  \textbf{AvP} is converging
linearly:

\begin{theorem}
\label{th:conv_GA} If  the linear regularity property
\eqref{linreg0} holds, then  we have  the following linear
convergence rate for the last iterate $x^k$ generated by algorithm
\textbf{AvP} with the optimal stepsize $\alpha_k = 1/\Uconst$:
\begin{equation}
\label{AvP_linconv} \dist^2_{\cX}(x^{k+1}) \leq \left(1-
\frac{1}{\Uconst \cdot \kappa } \right) \dist^2_{\cX}(x^k),
\end{equation}
or in terms of function values:
\begin{equation*}
F(x^{k}) - F^* \leq    \left(1- \frac{1}{\Uconst \cdot \kappa }
\right)^k \frac{\Uconst \dist^2_{\cX}(x^0)}{2}.
\end{equation*}
\end{theorem}

%\begin{proof}
\noindent \textit{Proof}:  Let $x^*$ be any element of $\cX$. Then, we have the following:
\begin{align}
\label{eq:decreaseGA} & \|x^{k+1} - x^*\|^2 \!  = \left\|x^k - x^*
- \alpha_k \left( x^k - \Exp{\Pi_{\cX_S}(x^k)} \right)\right\|^2 \nonumber \\
&=  \|x^k -x^*\|^2 - 2\alpha_k   \left\langle x^k - x^* , x^k -
\Exp{\Pi_{\cX_S}(x^k)} \right \rangle + \alpha_k^2
\left\| x^k - \Exp{\Pi_{\cX_S}(x^k)} \right\|^2 \nonumber\\
& \overset{\eqref{eq:gamma}}{\leq} \!\! \|x^k \!-x^*\|^2 \!-
2\alpha_k \!\left\langle x^k \!- x^*, x^k \!- \Exp{\Pi_{\cX_S}(x^k)}
\right \rangle +\! \Uconst \alpha_k^2
\Exp{ \left\| x^k \!- \Pi_{\cX_S}(x^k) \right\|^2} \nonumber\\
& = \|x^k \!-x^*\|^2 - 2 \alpha_k \Exp{\left\langle x^k -
\Pi_{\cX_S}(x^k) + \Pi_{\cX_S}(x^k) - x^*, x^k - \Pi_{\cX_S}(x^k)
\right\rangle} \nonumber \\
&  \qquad \qquad + \Uconst \alpha_k^2
\Exp{ \left\| x^k \!- \Pi_{\cX_S}(x^k) \right\|^2} \nonumber\\
& =  \|x^k \!-x^*\|^2 - 2 \alpha_k \Exp{ \left\| x^k \!-
\Pi_{\cX_S}(x^k) \right\|^2}  + \Uconst \alpha_k^2 \Exp{ \left\|
x^k \!- \Pi_{\cX_S}(x^k) \right\|^2} \nonumber \\
&  \qquad \qquad - 2 \alpha_k \Exp{\left\langle \Pi_{\cX_S}(x^k) -
x^*, x^k - \Pi_{\cX_S}(x^k) \right\rangle} \nonumber \\
& \leq  \|x^k \!-x^*\|^2 - (2 \alpha_k - \Uconst \alpha_k^2) \Exp{
\left\| x^k \!-\Pi_{\cX_S}(x^k) \right\|^2} \nonumber \\
& = \|x^k \!-x^*\|^2 - 2 (2 \alpha_k - \Uconst \alpha_k^2) F(x^k).
\end{align}
where the second inequality follows from the optimality condition of
the projection  $ \left\langle \Pi_{\cX_S}(x^k) - x^*, x^k -
\Pi_{\cX_S}(x^k) \right\rangle \geq 0$ for all $x^* \in \cX
\subseteq \cX_{S}$. From \eqref{eq:decreaseGA} we observe that the
fastest decrease is obtained from maximizing $2 \alpha_k - \Uconst
\alpha_k^2$, which leads to the optimal stepsize $\alpha_k =
1/\Uconst$. For this choice of the  stepsize,  $\alpha_k =
1/\Uconst$, and using $F^*=0$ we obtain from \eqref{eq:decreaseGA}:
\begin{align*}
\|x^{k+1} - x^*\|^2 & \leq \|x^k \!-x^*\|^2 - \frac{2}{\Uconst}
(F(x^k) - F^*) \\
& \overset{\eqref{qg_Lip3}}{\leq} \|x^k \!-x^*\|^2 -
\frac{1}{\Uconst \kappa} \|x^k \!-x^*\|^2,
\end{align*}
which implies immediately the statement of the theorem. \qed
%\end{proof}

\noindent Note that from the proof of Theorem \ref{th:conv_GA} it
follows that we can achieve linear convergence for the last iterate
generated by algorithm \textbf{AvP} with stepsizes satisfying $0 <
\alpha_k < 2/\Uconst$ Moreover, $\Uconst \kappa$ represents the
condition number of the convex feasibility problem
\eqref{convexfeas_original} or of its stochastic reformulation
\eqref{eq:reform_stoch_opt} (see Theorem \ref{th_shapeF}).

\noindent Let us consider finding a point in the finite intersection
of convex sets $(\cX_i)_{i \in [m]}$, that is  $\cX = \cap_{i=1}^m
\cX_i$. Further, we consider a uniform probability  on $\Omega=[m]$
and we choose the minibach sample size $N=m$, then the average
projection algorithm \textbf{AvP} becomes the barycentric method:
\[\textbf{AvP}(1/m): \quad  x^{k+1} = x^k - \alpha_k \left( x^k -
\frac{1}{m} \sum_{i=1}^m \Pi_{\cX_i}(x^k) \right).  \]
The barycentric method  was shown  to converge asymptotically to a
point in the intersection of the closed  convex sets $(\cX_i)_{i \in
[m]}$, see e.g. \cite{Com:97}. Recall that we denoted  $D =
\text{diag}(\|A_1\|^{-2}, \cdots, \|A_m\|^{-2})$. Let us derive
convergence rates for the  barycentric method \textbf{AvP}$(1/m)$
for two particular cases of sets:

\noindent \textbf{(i)}: Consider the problem of finding a solution
to a linear system $Ax = b$, where $A$ is an $m \times n$ matrix. In
this case $\cX_i = \{x: \ A_i^T x = b_i\}$. Then, from Theorem
\ref{th:conv_GA} the barycentric method \textbf{AvP}$(1/m)$ with the
optimal stepsize $\alpha_k = 1/\Uconst$  converges linearly:
\begin{align*}
\Exp{\dist_{\cX}^2(x^{k})} & \overset{\eqref{AvP_linconv}}{\leq}
\left( 1 - \frac{1}{\Uconst \kappa} \right)^k \dist_{\cX}^2(x^0) \\
& \overset{\eqref{cor1_Ules2}+\eqref{const_lineq}}{=} \left( 1 -
\frac{\lambda_{\min}^{\text{nz}}\left(A^T D
A\right)}{\lambda_{\max}(A^T D A)} \right)^k \dist_{\cX}^2(x^0).
\end{align*}

\noindent \textbf{(ii)}: Consider now the more general problem of
finding a solution to a  system of linear inequalities $Ax \leq b$,
where $A$ is an $m \times n$ matrix. Then $\cX_i = \{x: \ A_i^T x
\leq b_i\}$. From Theorem \ref{th:conv_GA}  it follows that the
barycentric method \textbf{AvP}$(1/m)$ with the optimal stepsize
$\alpha_k = 1/\Uconst$ converges also linearly:
\begin{align*}
\Exp{\dist_{\cX}^2(x^{k})} & \overset{\eqref{AvP_linconv}}{\leq}
\left( 1 - \frac{1}{\Uconst \kappa} \right)^k \dist_{\cX}^2(x^0) \\
& \overset{\eqref{cor2_Ules2}+\eqref{const_linineq_u}}{=} \left( 1 -
\frac{1}{\max\limits_{i=1:m} \|A_i\|^2 \lambda_{\max}(A^T D A)
\tilde \kappa} \right)^k \dist_{\cX}^2(x^0).
\end{align*}
Note that from Theorem \ref{th:conv_GA} it follows immediately that
the basic barycentric method $x^{k+1} = \frac{1}{m} \sum_{i=1}^m
\Pi_{\cX_i}(x^k)$, i.e. stepsize $\alpha_k = 1$, converges linearly:
\begin{align*}
\Exp{\dist_{\cX}^2(x^{k})} & \overset{\eqref{AvP_linconv}}{\leq}
\left( 1 - \frac{2-\Uconst}{\kappa} \right)^k \dist_{\cX}^2(x^0).
\end{align*}

\vspace{5pt}

\noindent However, our algorithmic framework leads to new schemes.
For example, for a general  probability distribution $(p_i)_{i \in
[m]}$ on $\Omega=[m]$ and $N=m$, the average projection algorithm
\textbf{AvP} has the iteration:
\[\textbf{AvP}(p_i): \quad  x^{k+1} = x^k - \alpha_k \left( x^k -
\sum_{i=1}^m p_i \Pi_{\cX_i}(x^k) \right).  \]
If we choose the probabilities $p_i =
\frac{\norm{A_i}^2}{\norm{A}^2_F}$, then this method has the
following convergence rates  for  linear systems and linear
inequalities:

\noindent \textbf{(iii)}: For a linear system $Ax = b$,  from
Theorem \ref{th:conv_GA} the
\textbf{AvP}$(\norm{A_i}^2/\norm{A}^2_F)$ method with the optimal
stepsize $\alpha_k = 1/\Uconst$  converges linearly:
\begin{align*}
\Exp{\dist_{\cX}^2(x^{k})} & \overset{\eqref{AvP_linconv}}{\leq}
\left( 1 - \frac{1}{\Uconst \kappa} \right)^k \dist_{\cX}^2(x^0) \\
& \overset{\eqref{cor1_Ules2}+\eqref{const_lineq}}{=} \left( 1 -
\frac{\lambda_{\min}^{\text{nz}}\left(A^T
A\right)}{\lambda_{\max}(A^T  A)} \right)^k \dist_{\cX}^2(x^0).
\end{align*}

\noindent \textbf{(iv)}: For a  system of linear inequalities $Ax
\leq b$, from Theorem \ref{th:conv_GA}  the previous method  with
the optimal stepsize $\alpha_k = 1/\Uconst$ converges also linearly:
\begin{align*}
\Exp{\dist_{\cX}^2(x^{k})} & \overset{\eqref{AvP_linconv}}{\leq}
\left( 1 - \frac{1}{\Uconst \kappa} \right)^k \dist_{\cX}^2(x^0) \\
& \overset{\eqref{cor2_Ules2}+\eqref{const_linineq}}{=} \left( 1 -
\frac{1}{ \lambda_{\max}(A^T  A) \tilde \kappa} \right)^k
\dist_{\cX}^2(x^0).
\end{align*}

% \noindent However, for general  convex sets $(\cX_i)_{i \in [m]}$ it
% is difficult to accurately estimate the Lipschitz constant
% $\Uconst$. It that case we can estimate  adaptively the Lipschitz
% constant:
% \[ \Uconst_k = \frac{\| \sum\limits_{i=1}^{m}  p_i(x^ k -
% \Pi_{\cX_{i}}(x^k))\|^2}{ \sum\limits_{i=1}^{m} p_i \left\| x^k -
% \Pi_{\cX_{i}}(x^k) \right\|^2 } \leq \Uconst \] and use an adaptive
% stepsize $\alpha_k = \alpha/\Uconst_k$ for some $\alpha \in (0, \
% 2)$. Then, we recover an adaptive variant of \textbf{AvP} analyzed
% in \cite{Com:97} and whose effectiveness has been shown in
% \cite{CenChe:12}. Note however, that in \cite{Com:97} only
% asymptotic convergence results are provided, and no rates of
% convergence, while with our analysis we can provide explicit rates.

%%%%%%%%%%%%%%%%%%%%%%%%%%%%%%%%%%%%%%%%%%%%%%%%%%%%%%%%%%%%%%%%%%%%%%%%%%%

\subsection{Stochastic Alternating Projection algorithm: $N=1$}
In this section we analyze in more detail a particular case of
scheme \textbf{SPA} which uses a single projection for the updates.
That is, in  \textbf{SPA} we choose $N=1$, which results in the
Stochastic Alternating Projection (\textbf{SAP}) scheme:

\vspace{10pt}

\begin{center}
\framebox{
\parbox{8.5 cm}{
\begin{center}
\textbf{ Algorithm SAP  }
\end{center}
Choose $x^0 \in \R^n$ and positive stepsizes $\{\alpha_k\}_{k \ge 0}$\\
For $k\geq 0$ repeat:
\begin{enumerate}
\item Choose randomly a sample $S_k \sim \Prob$
\item Compute $x^{k+1} = x^k - \alpha_k \left(x^k -   \Pi_{\cX_{S_k}}(x^k)\right)$
\end{enumerate}
}}
\end{center}

\vspace{10pt}

\noindent Algorithm \textbf{SAP} can be viewed as a random
implementation of the alternating projection method, which generates
a sequence of iterates by projecting on the sets cyclically. The
alternating projection algorithm has been proposed by Von Neumann
\cite{Neu:50} for the intersection problem of two subspaces in a
Hilbert space, and it has many generalization and extensions
\cite{BauNol:13,DeuHun:06,Ned:10}. A nice survey of the work in this
area is given in \cite{BauBor:96}. The first convergence rate result
for the alternating projection algorithm under the assumption that
the  intersection set has a nonempty interior has been given in
\cite{GubPol:66}. Unlike the alternate projection method (which is
deterministic), the  algorithm \textbf{SAP}  utilize random
projections. The convergence rate of \textbf{SAP} for a finite
intersection of simple convex sets has been given recently in
\cite{Ned:10,Ned:11}. From the convergence analysis of previous
section it follows that the stepsize in \textbf{SAP} can be chosen
as:
\[  \delta \leq \alpha_k \leq 2 -\delta, \]
since for $N=1$ we have $\Uconst_N=1$. Moreover,  the optimal
stepsize is $\alpha_k=1$. However, it has been observed in practice
that overrelaxations, that is $\alpha_k \in [1, \ 2]$, make
\textbf{SAP} to perform better. Further note that for specific sets
and probabilities we recover well known algorithms from
literature:\\

\noindent \textbf{(i)}: Consider the  problem of finding a solution
to a linear system $Ax = b$, where $A$ is an $m \times n$ matrix.
Further, assume $\Omega=\{e_1, \cdots, e_m\}$ and the probability
distribution $\Prob(S=e_i) = \frac{\norm{A_i}^2}{\norm{A}^2_F}$.
Then, \textbf{SAP} with $\alpha_k=1$ is the randomized Kaczmarz
algorithm from \cite{StrVer:09}:
\[ x^{k+1} = x^k -  \frac{A_i^T x^k-b_i}{\|A_i\|^2} A_i.  \]
Moreover, for these choices of the probabilities and stepsize, our
convergence analysis matches exactly the one in \cite{StrVer:09},
that is \textbf{SAP} is converging linearly:
\begin{equation*}
\Exp{\dist_{\cX}^2(x^{k})} \overset{\eqref{SPA_linconv}}{\leq}
\left( 1 - \frac{1}{\kappa} \right)^k \dist_{\cX}^2(x^0)
\overset{\eqref{const_lineq}}{=} \left( 1 -
\frac{\lambda_{\min}^{\text{nz}}\left(A^TA\right)}{\norm{A}^2_F}
\right)^k \dist_{\cX}^2(x^0).
\end{equation*}

\noindent However, \textbf{SAP} generalizes the randomized Kaczmarz
algorithm from \cite{StrVer:09}, considering for a random matrix
$S_k \in \rset^{m \times q}$ the general iteration:
\[  x^{k+1} = x^k - \alpha_k  A^TS_k(S^T_kAA^TS_k)^{\dagger} S^T_k(Ax^k-b).
\]
Notice that for constant stepsize $\alpha_k = 1$, the previous
\textbf{SAP} scheme is equivalent with  the randomized iterative
method of \cite{GowRic:15}. For this choice of the stepsize, our
convergence analysis matches exactly the one in \cite{GowRic:15}:
\begin{align*}
\Exp{\dist_{\cX}^2(x^{k})}  & \overset{\eqref{SPA_linconv}}{\leq}
\left( 1 - \frac{1}{\kappa} \right)^k \dist_{\cX}^2(x^0) \\
& \overset{\eqref{kappa_linear}}{=} \left(1 -
\lambda_{\min}^{\text{nz}}(A^T \Exp{ S(S^TAA^TS)^{\dagger} S^T} A)
\right)^k \dist_{\cX}^2(x^0).
\end{align*}

\vspace{5pt}

\noindent \textbf{(ii)}: Consider now the more general problem of
finding a solution to a  system of linear inequalities $Ax \leq b$,
where $A$ is an $m \times n$ matrix. Further, assume as above
$\Omega=\{e_1, \cdots, e_m\}$ and the probability distribution
$\Prob(S=e_i) = \frac{\norm{A_i}^2}{\norm{A}^2_F}$. Then,
\textbf{SAP} with $\alpha_k=1$ is the  Algorithm 4.6 from
\cite{LevLew:10}:
\[ x^{k+1} = x^k -  \frac{\Pi_+(A_i^T x^k-b_i)}{\|A_i\|^2} A_i.  \]
For  these choices of the probabilities and stepsize, our
convergence analysis matches exactly the one in \cite{LevLew:10}:
\begin{equation*}
\Exp{\dist_{\cX}^2(x^{k})} \overset{\eqref{SPA_linconv}}{\leq}
\left( 1 - \frac{1}{\kappa} \right)^k \dist_{\cX}^2(x^0)
\overset{\eqref{const_linineq}}{=} \left( 1 -
\frac{1}{\tilde{\kappa} \norm{A}^2_F} \right)^k \dist_{\cX}^2(x^0).
\end{equation*}

\noindent However, \textbf{SAP} generalizes the Algorithm 4.6 from
\cite{LevLew:10}, considering for a random vector $S_k \in
\rset^{m}_+$ the general  iteration:
\[  x^{k+1} = x^k - \alpha_k  \frac{\Pi_+(S^T_kAx -
S^T_kb)}{\norm{A^TS_k}^2}A^TS_k.
\]
Under the settings of Theorem \ref{lemma_halfspace} we obtain:
\begin{align*}
\Exp{\dist_{\cX}^2(x^{k})}  & \overset{\eqref{SPA_linconv}}{\leq}
\left( 1 - \frac{1}{\kappa} \right)^k \dist_{\cX}^2(x^0)
\overset{\eqref{const_linineq}}{=} \left(1 - \frac{1}{\tilde{\kappa}
\sum_{S \in \Omega_r} \norm{A^T S}^2} \right)^k \dist_{\cX}^2(x^0).
\end{align*}

\vspace{5pt}

\noindent \textbf{(iii)}: Finally, we can consider the convex
feasibility problem where the intersection set has a nonempty
interior. First, let us investigate when the sequence $\|x^k -
x^*\|$ is decreasing. For $N=1$ and $\alpha_k \in [0, \ 2]$ it
follows from \eqref{eq:decrease} that:
\[ \|x^{k+1} - x^*\|^2 \leq \|x^k - x^*\|^2 - (2 \alpha_k - \alpha_k^2)
\| x^k - \Pi_{\cX_{S^k}}(x^k)\| \quad \forall k \geq 0, \] that is
the sequence $\norm{x^k - x^*}$ is nonincreasing. Similarly, for $N
\geq 1$ and $\alpha_k \in [0, \ 1]$ it follows that:
\begin{align*}
\|x^{k+1} - x^*\|^2 & = \left\| (1 -\alpha_k) (x^k - x^*) + \alpha_k
\left(
\frac{1}{N}\sum\limits_{i=1}^N \Pi_{\cX_{S^k_i}}(x^k) -x^* \right) \right\|^2 \\
& \leq (1 - \alpha_k) \|x^k - x^* \|^2 + \alpha_k \left\|
\frac{1}{N}\sum\limits_{i=1}^N \Pi_{\cX_{S^k_i}}(x^k) -x^*
\right\|^2\\
& \leq (1 - \alpha_k) \|x^k - x^* \|^2 + \frac{\alpha_k}{N}
\sum\limits_{i=1}^N \|\Pi_{\cX_{S^k_i}}(x^k) -x^*\|^2 \\
& = \|x^k - x^*\|^2  + \frac{\alpha_k}{N} \sum\limits_{i=1}^N
\left(\|\Pi_{\cX_{S^k_i}}(x^k) -x^*\|^2 - \|x^k - x^* \|^2 \right)\\
& \leq \|x^k - x^*\|^2 - \frac{\alpha_k}{N} \sum\limits_{i=1}^N
\|x^k - \Pi_{\cX_{S^k_i}}(x^k)\|^2 \quad \forall k \geq 0.
\end{align*}
The last inequality follows from the bound $\| x^k -
\Pi_{\cX_{S^k_i}}(x^k) \|^2 + \|\Pi_{\cX_{S^k_i}}(x^k) - x^*\|^2
\leq \|x^k - x^* \|^2 $ for all $x^* \in \cX \subseteq \cX_{S_i^k}$.
Therefore, for $N \geq 1$ and $\alpha_k \in [0, \ 1]$ we also have a
nonincreasing  sequence $\norm{x^k - x^*}$. In conclusion, for the
two choices for $N$ and $\alpha_k$ given above we have:
\begin{equation*}
\norm{x^{k} - x^*} \le  \norm{x^0 - x^*} \quad \forall x^* \in \cX,
\; k \ge 0.
\end{equation*}
An important application of the previous inequality is that when the
set $\cX$ contains a ball with radius $\delta$ centered in
$\bar{x}$. By taking $x^* = \bar{x}$ in the previous relation, we
have: $\norm{x^{k} - \bar{x}} \le \norm{x^0 - \bar{x}}$ for all $k
\ge 0$. This implies that under the settings of Theorem
\ref{lemma_ball}, one should choose the compact set $Q=\{x: \;
\norm{x - \bar{x}} \le \norm{x^0 - \bar{x}} \}$, such that the
linear regularity constant given  in \eqref{const_ball} becomes:
\begin{equation}
\label{const_ball1} \kappa = \frac{ \norm{x^0 - \bar{x}}^2}{\delta^2
\min\limits_{S \in \Omega} \; p_S},
\end{equation}
since all the points of interest for which the linear regularity
property has to hold are the iterates $\{x^k\}_{k \ge 0}$. Then,
\textbf{SAP} with $\alpha_k=1$ is the random projection algorithm
from \cite{Ned:10}. For this choice of the stepsize and under the
setting of Theorem  \ref{lemma_ball}, the algorithm \textbf{SAP}
attains the following linear rate:
\begin{equation*}
\Exp{ \dist_{\cX}^2(x^{k}) }  \overset{\eqref{SPA_linconv}}{\leq}
\left( 1 - \frac{1}{\kappa} \right)^k \dist_{\cX}^2(x^0)
\overset{\eqref{const_ball}+\eqref{const_ball1}}{=} \left(1 -
\frac{p_{\min}\delta^2 }{R^2} \right)^k \dist_{\cX}^2(x^0),
\end{equation*}
where $p_{\min} = \min_{S \in \Omega} p_S$ and  $R = \norm{x^0 -
\bar{x}}$. A similar convergence rate has been derived in
\cite{Ned:10} for this particular scheme.

\section{Conclusions}
We have proposed new stochastic reformulations of the classical convex feasibility problem and analyzed the problem conditioning parameters in relation with (linear) regularity assumptions on the individual convex sets. Then, we have introduced a  general random projection algorithmic framework, which extends to the random
settings many existing projection schemes, designed for the general convex feasibility problem. Based on the conditioning parameters, besides the asymptotic convergence results, we have also derived explicit sublinear and linear convergence rates for this general algorithm. The convergence rates show specific dependence on the number of projections averaged at each iteration. Our general random projection algorithm also allows to project simultaneously on several sets, thus providing great flexibility in matching the implementation of the algorithms on the parallel architecture at hand.

%%%%%%%%%%%%%%%%%%%%%%%%%%%%%%%%%%%%%%%%%%%%%%%%%%%%%%%%%%%%%%%%%%%%%%%%%
%%%%%%%%%%%%%%%%%%%%%%%%%%%%%%%%%%%%%%%%%%%%%%%%%%%%%%%%%%%%%%%%%%%%%%%%%%

\end{document}